\documentclass[10pt]{amsart}

\usepackage{amsthm,amsmath,amssymb,amscd}
\usepackage{amssymb}
\usepackage{graphics}

\newtheorem{theorem}{Theorem}[section]

\newtheorem{remark}[theorem]{Remark}
\newtheorem{lemma}[theorem]{Lemma}
\newtheorem{proposition}[theorem]{Proposition}
\newtheorem{definition}{Definition}[section]

\newcommand{\R}{\mathbb{R}}

\begin{document}

\title [A variational characterization of $J$-holomorphic curves]
 {A variational characterization of $J$-holomorphic curves}
\author{Claudio Arezzo, Jun Sun}

\address{The Abdus Salam International Centre for Theoretical Physics\\ Trieste, Italy
and \newline Universit\'a di Parma, Italy.}
\email{arezzo@ictp.it}

\address{CIRM, Fondazione Bruno Kessler, Via Sommarive, 14 - Povo, I-38123 \\ Trento (TN), Italy.}
\email{sunjun@fbk.eu}


\date{}

\maketitle
\begin{abstract}
In this paper, we prove that if the area functional of a surface $\Sigma^2$ in a symplectic manifold $(M^{2n},\bar{\omega})$ has a critical point or has a compatible stable point in the same cohomology class, then it must be $J$-holomorphic. Inspired by a classical result of Lawson-Simons, we show how various restrictions of the stability assumption to variations of metrics in the space "projectively induced" metrics are enough to  give the desired conclusion. \end{abstract}

\vspace{.2in}

\section{Introduction}

\allowdisplaybreaks

\vspace{.1in}

\noindent By a well-known extension of Wirtinger's inequality we know that complex submanifolds of K\"ahler manifolds minimize volume in their homology class. A classical problem is to determine to which extent the converse holds.
For example, Lawson-Simons (\cite{LS}) proved that any stable minimal surface in $\textbf{CP}^n$ is holomorphic. Siu-Yau proved the same result when the ambient manifold has a metric of positive holomorphic bisectional curvature (\cite{SY}). Micallef (\cite{Mi}) studied complete stable minimal surfaces in $\textbf{R}^4$, and proved that, under some further assumptions, stable minimal surfaces must be holomorphic with respect to some complex structure on $\textbf{R}^4$. He also proved the analogue result for compact surfaces in flat $4$-tori (\cite{M2}).

\noindent Since then a series of examples of stable symplectic minimal surfaces (representing $(1,1)$-classes in homology) which are {\em{not}} holomorphic w.r.t. any complex structure have been found under different curvature assumptions of the ambient metric (see e.g. \cite{AM} for higher dimensional flat tori,
\cite{AMP} for higher dimensional euclidean spaces, \cite{Ar} for K\"ahler-Einstein surfaces of negative curvature, \cite{MW} for K3 surfaces and \cite{AL} for
K\"ahler-Einstein manifolds of dimensions greater than $4$ and positive curvature).

\noindent In this paper we consider immersions of surfaces into symplectic manifolds endowed with a compatible almost complex structure which are critical and stable w.r.t. variations of the ambient metric. This seems a very natural extension of the classical situation since we are using the metric just as a tool to detect $J$-holomorphicity of a submanifold but we are not really interested in any of its  riemannian properties.

\noindent Allowing arbitrary deformations of the metric on the ambient space give rise to a badly defined concept of critical point, as explained in the Appendix A. Instead, as we have a symplectic form $\bar{\omega}$ on $M$, we restrict ourselves to variations of the metric induced via the (tamed almost) complex structure $J$ by special variations of the symplectic form in the same cohomology class as $\bar{\omega}$.

\vspace{.1in}

\noindent Let us first recall some notations. For a compact symplectic manifold $(M^{2n},\bar\omega)$, it is known that (Corollary 12.7 of \cite{Ca}) there always exists an almost complex structure $J$ that is tamed and compatible with $\bar\omega$. Namely, $\bar\omega(X,JX)>0$ for $0\neq X\in TM$ and $\bar{\omega}(JX,JY)=\bar{\omega}(X,Y)$. Thus we can define the associated Riemannian metric by
\begin{equation}\label{metric-compatible}
    \bar{g}(X,Y)=\bar{\omega}(X, JY).
\end{equation}
Given the almost complex structure $J$, we have the splitting of the tangent space $TM^{\textbf{C}}=T^{1,0}M\oplus T^{0,1}M$ and the splitting of the cotangent space $\Lambda^1_{\textbf{C}}M=\Lambda^{1,0}M\oplus \Lambda^{0,1}M$. Then by definition, $\partial\psi$ and $\bar{\partial}\psi$ are just the components of $d\psi\in \Lambda^1M$ in $\Lambda^{1,0}M$ and $\Lambda^{0,1}M$, respectively. Set
\begin{equation}\label{E2.2}
    d=\partial+\bar{\partial}, \ \ d^c=\sqrt{-1}(\bar{\partial}-\partial).
\end{equation}

\noindent Let $\mathcal{H} = \{  \rho \in \mathcal{C}^{\infty}(M,\R) \mid  \bar\omega_{\rho}:=\bar\omega+dd^c\rho\,\,  \mbox{tames} \,\, J\}$, which is clearly a nonempty open subset of $\mathcal{C}^{\infty}(M,\R)$. 

\noindent To each $\rho \in \mathcal{H}$ we can associate a riemannian metric $\bar{g}_{\rho}$ on $M$ defined by

\begin{equation}\label{e2.1}
 \bar{g}_{\rho}(X,Y)=\frac{1}{2}\left(\bar{\omega}_{\rho}(X, JY)+\bar{\omega}_{\rho}(Y, JX)\right)
\end{equation}

\vspace{.1in}

\noindent Let $\Sigma$ be a closed real surface and
\begin{equation*}
    F:\Sigma \to M
\end{equation*}
be an immersion. 
We can then define 

\begin{equation}
\label{defarea}
\mathcal{A}(\rho) = \mbox{Area}(F(\Sigma), F^*(\bar{g}_{\rho})) = \int_{\Sigma}d\mu_{\rho} \,\, ,
\end{equation}

where $d\mu_{\rho}$ is the volume form of the induced metric $g_{\rho}:=  F^*(\bar{g}_{\rho})$.

\noindent Clearly the functional $\mathcal{A}$ depends only on the metric $\bar{g}_{\rho}$ (which, by (\ref{e2.1}) depends in turn on 
$\bar{\omega}_{\rho}$) and not on the choice of refence metric $\bar{\omega}$ and the potential $\rho$. For this reason we will often think of $\mathcal{A}$ as a functional on the ``tamed" subset of $[\bar{\omega}]$.

\vspace{.1in}

\begin{definition}\label{stationary}
Given an immersion $F:\Sigma^2\to (M,\bar\omega,J,\bar g)$, we say that  the area functional $\mathcal{A}$ has a critical point $\rho \in \mathcal{H}$ if for any $\phi(t) \in \mathcal{H}$ with $\phi(0) = \rho$\begin{equation*}
    \mathcal{A}'(0)=0.
\end{equation*}
\end{definition}

\noindent It is a simple consequence of Stokes' Theorem that if $\Sigma$ is $J$-holomorphic (even in the more general {\em{tamed}} situation), the functional $\mathcal{A}$ is constant on $\mathcal{H}$.
The first result in this paper shows that the existence of a critical point of $\mathcal{A}$ is enough to guarantee the $J$-holomorphicity:

\vspace{.1in}

\noindent \textbf{Theorem 2.5:} {\it Let $(M^{2n},\bar{\omega},J)$ be a compact symplectic manifold with compatible almost complex structure $J$ and $F:\Sigma^2\to M$ be an immersion. If the area functional ${\mathcal A}$ has a critical point in $\mathcal{H}$, then the immersion is $J$-holomorphic.}

\vspace{.1in}

\noindent In light of our knowledge about the relationship between stable minimal surfaces and holomorphic curves, it is natural to look at special properties of the second variation of the functional $\mathcal{A}$:

\begin{definition}\label{A-stable}
Given an immersion $F:\Sigma^2\to (M,\bar\omega,J,\bar g)$, we say that $\rho \in \mathcal{H}$ is a stable point for  the area functional ${\mathcal A}$ if 
\begin{equation*}
    {\mathcal A}''(0)\geq 0
\end{equation*}
for any $\phi(t)\in \mathcal{H}$, $\phi(0)=\rho$. Furthermore, if $J$ is compatible with $\bar\omega_{\rho}$, then we say $\rho$ is a  \textbf{compatible stable point}.
\end{definition}

\noindent Note that the definition of ${\mathcal A}$-stability (as well as all the other stability notions we are going to study) does not require $\rho$  to be a critical point of the area functional. Our next result shows that the existence of a compatible stable point is also enough to guarantee the $J$-holomorphicity:

\vspace{.1in}

\noindent \textbf{Theorem 3.2:} {\it Let $(M^{2n},\bar{\omega},J)$ be a compact symplectic manifold with compatible almost complex structure $J$ and $F:\Sigma^2\to M$ be an immersion. If the area functional ${\mathcal A}$ has a compatible stable point $\rho \in \mathcal{H}$, then the immersion is $J$-holomorphic.}

\vspace{.1in}

\noindent As above, the converse is also true even without assuming that $J$ is compatible with $\bar\omega$.

\vspace{.1in}

\noindent Checking the proof of the above theorem carefully, we see that the result is also true in the complete noncompact case. In particular, this applies to minimal submanifold in $\textbf{R}^{2n}$. In this case, we need the test function to have compact support.

\vspace{.1in}

\noindent \textbf{Theorem 3.4:} {\it Let $(M^{2n},\bar{\omega},J)$ be a complete noncompact symplectic manifold with compatible almost complex structure $J$ and $F:\Sigma^2\to M$ be an immersion. If the area functional ${\mathcal A}$ has a compatible stable point $\rho \in \mathcal{H}$, then the immersion is $J$-holomorphic.}

\vspace{.1in}

\noindent The above Theorems show an intriguing analogy with of a classical picture discovered by Sacks-Uhlenbeck \cite{SU} and Schoen-Yau \cite{ScY} to construct minimal surfaces, which, as they proved, can be generated by first fixing a {\em{ metric on the domain surface}} and finding an energy minimizing map, and then letting the metric on the base vary. In our case minimizing w.r.t. the {\em{metric on the target}}  plays the corresponding role which not surprisingly encodes a similar strategy since of course for what the area can detect the immersion is not an independent parameter compared to the ambient metric.

\vspace{.1in}

\noindent It is clear that, while very simple to state, the induced existence problem is very difficult to attack, since we introduced a parameter (the target metric) which varies freely in an infinite dimensional noncompact space (the ``K\"ahler potentials" $\mathcal{H}$). We then try to identify some geometrically meaningful finite dimensional subspaces or submanifolds of metrics which are enough to detect holomorphicity.

\noindent When $M$ is an algebraic manifold we can embed it into some complex projective space ${\textbf{CP}^N}$ holomorphically. Lawson-Simons' result (\cite{LS}) tells us that a submanifold is holomorphic if the second variation of the area functional (built with the metric induced by the projective space) is nonnegative under the holomorphic deformation of $M$ in ${\textbf{CP}^N}$. The latter means that the target metric varies in some finitely dimensional submanifold of metrics representing the original K\"aher class. Their result is true for submanifolds of any even  dimension but only for projectively induced K\"ahler metrics (in particular representing a rational class).

\vspace{.1in}

\noindent We first consider the same case as Lawson-Simons' under a slightly different stability assumption. We call the area functional has a \emph{linearly projectively stable point} if the variation of the metric on the target is linear in time along the directions induced by holomorphic deformations of the complex projective space (Lawson-Simons' assumption can be stated as to require the variations of the metric to live in this subspace for all time).

\vspace{.1in}
\noindent \textbf{Theorem 4.1:} {\it Let $(M,\bar\omega,J,\bar g)$ be an algebraic manifold with all structures induced by the projective space as above and $F:\Sigma^2\to M$ be an immersion. If the area functional has a linearly projectively stable point, then the immersion must be holomorphic with respect to the complex structure $J$.}

\vspace{.1in}

\noindent This can be seen as a mild modification of Lawson-Simons' result in the surface case, yet our proof differs significantly from theirs in that we explicitly identify in the nonholomorphic case a Killing field which induces an area-decreasing variation, while they had to work on the whole space of Killing fields and using heavily the homogeneous space structure of the projective space to average the variations of area.

\vspace{.1in}

\noindent The advantage of our proof of Theorem $5.1$ is that it generalizes to any symplectic manifold with rational symplectic class and to any complex projective manifold with any K\"ahler form. Indeed, let $(M,\bar\omega,J)$ be a symplectic manifold with rational symplectic class and compatible almost complex structure. It is known that (\cite{BU}), there exists an approximately $J$-holomorphic embedding of $M$ in to some complex projective space $\textbf{CP}^{N_k}$. In the symplectic case, using as above the holomorphic deformations of $\textbf{CP}^{N_k}$, we can extend the notion of linearly projectively stable point to that of compatible linearly ${\mathcal A}^k$-stable point, which again means that the target symplectic form varies (still linearly) in the direction of projectively induced forms.

\vspace{.1in}
\noindent \textbf{Theorem 5.1:} {\it
Let $(M^{2n},\bar{\omega},J_M,\bar g )$ be a symplectic manifold as above and $F:\Sigma^2\to M$ be an immersion. There exists an integer $K_1$, such that if the area functional has a compatible linearly ${\mathcal A}^k$-stable point for some $k\geq K_1$, then the immersion must be $J_M$-holomorphic.}

\vspace{.1in}

\noindent Using our second variation formula, we can show that for surface case, we can obtain similar result for algebraic manifolds but whose K\"ahler form represents any real class. We can define the notion of {\em{$k$-linearly projectively stable}} point, which means that the second variation of the area functional is nonnegative when the potential of the target metric varies along the directions in the finitely dimensional linear subspace of the space of potentials determined by the Killing vector fields of ${\textbf{CP}^{N_k}}$, where this projective space is the ambient of a diagonal approximating sequence of rational forms converging to the original class and the above approximation described in the symplectic case (which in this case relies on a famous Theorem by Tian \cite{Tian}).

\vspace{.1in}

\noindent \textbf{Theorem 6.1:} {\it Let $(M,J)$ be an algebraic manifold, and $\bar\omega$ be any K\"ahler metric with $[\bar\omega]\in H^2(M,\textbf{R})\cap H^{1,1}(M,\textbf{C})$ and $F:\Sigma^2\to M$ be an immersion. There exists an integer $K_2$, such that if if the area functional has a $k$-linearly projectively stable point for some $k\geq K_2$, then the immersion must be holomorphic with respect to the complex structure $J$.}

\vspace{.1in}

\noindent We finally underline that most of our arguments and results are likely to carry over to higher dimensional submanifolds, and in fact to more general, less regular, objects. This will be clarified in a forthcoming paper.

\vspace{.1in}

\noindent The following sections are organized as follows: in Section 2 and Section 3 we compute the first and second variation formulas for the area functional under deformation of target metrics and prove the first two results; in Section 4, we study the case of an ambient algebraic manifold with induced structures and linearly projectively stable point; in Section 5 and 6, we prove the symplectic case with rational classes and the K\"ahler case with any real K\"ahler class.

\vspace{.2in}

\section{Critical points of $\mathcal{A}$ and $J$-holomorphicity}

\vspace{.1in}

\noindent With the setup given in the introduction, we now compute the first variation of the area functional ${\mathcal A}$ and prove the first theorem.

\vspace{.1in}

\noindent  Let $\{x_1,x_2\}$ be local coordinates on $\Sigma$ and $\bar{g}_{\rho}(t)$ a variation of  $\bar{g}_{\rho}$ coming from a $1$-parameter deformation of $\rho$ in $\mathcal{H}$. . Then
\begin{equation}\label{metric}
    g_{\rho,ij}(t)=\bar{g}_{\rho}(t)\left(\frac{\partial F}{\partial x_i},\frac{\partial F}{\partial x_j}\right)
    =\frac{1}{2}\left\{\bar{\omega}_{\rho}(t)\left(\frac{\partial F}{\partial x_i},J\frac{\partial F}{\partial x_j}\right)
              +\bar{\omega}_{\rho}(t)\left(\frac{\partial F}{\partial x_j},J\frac{\partial F}{\partial x_i}\right)\right\}.
\end{equation}
Set
\begin{equation}\label{nu}
    \nu_{\rho}(t)=\frac{\sqrt{det(g_{\rho,ij}(t))}}{\sqrt{det(g_{\rho,ij}(0))}}.
\end{equation}
Then $\nu_{\rho}(t)$ is well-defined independent of the choice of coordinate system. Furthermore,
\begin{equation}\label{e2.6}
    {\mathcal A}(t)=\int_{\Sigma}\sqrt{det(g_{\rho,ij}(t))}=\int_{\Sigma}\nu_{\rho}(t)\sqrt{det(g_{\rho,ij}(0))},
\end{equation}
and therefore
\begin{equation}\label{e2.7}
    \frac{d}{dt}|_{t=0}{\mathcal A}(t)=\int_{\Sigma}\frac{d}{dt}|_{t=0}\nu_{\rho}(t)\sqrt{det(g_{\rho,ij}(0))}.
\end{equation}
Denote $(g_{\rho}^{ij})=(g_{\rho,ij})^{-1}$. By (\ref{nu}), we have in any local coordinate $\{x_1,x_2\}$
\begin{eqnarray}\label{e2.9}
  \frac{d}{dt}|_{t=0}\nu_{\rho}(t)
   &=& \frac{1}{2}\sum_{i,j=1}^{2}g_{\rho}^{ij}(0)g_{\rho,ij}'(0) \nonumber\\
   &=& \frac{1}{2}\sum_{i,j=1}^{2}g_{\rho}^{ij}(0)\frac{d}{dt}|_{t=0}
         \left\{\bar{g}_{\rho}(t)\left(\frac{\partial F}{\partial x_i},\frac{\partial F}{\partial x_j}\right)\right\}\nonumber \\
   &=& \frac{1}{2}\sum_{i,j=1}^{2}g_{\rho}^{ij}(0)\frac{d}{dt}|_{t=0}\left\{\bar{\omega}_{\rho}(t)
         \left(\frac{\partial F}{\partial x_i},J\frac{\partial F}{\partial x_j}\right)\right\}\nonumber \\
   &=& \frac{1}{2}\sum_{i,j=1}^{2}g_{\rho}^{ij}(0)\bar\omega_{\rho}'(0)\left(F_{x_i},J F_{x_j}\right).
\end{eqnarray}
Therefore the first variation formula is given by
\begin{equation}\label{first}
    {\mathcal A}(0)=\frac{1}{2}\sum_{i,j=1}^{2}\int_{\Sigma}g_{\rho}^{ij}(0)\bar\omega_{\rho}'(0)\left(F_{x_i},J F_{x_j}\right)d\mu_{\rho}.
\end{equation}

\vspace{.1in}

\noindent For our later use, let's recall the following simple facts:

\begin{lemma}
\begin{enumerate}
\item
For any smooth function $\psi$ on $M$, we have
\begin{equation}\label{e2.12}
    d^c\psi=-d\psi\circ J.
\end{equation}
\item
For any smooth function $\psi$ on $M$ and any tangent vector fields $X,Y$ on $M$, we have
\begin{equation}\label{e2.13}
    (dd^c\psi)(X,Y)=-(\overline{\nabla}^2\psi)(X,JY)+(\overline{\nabla}^2\psi)(Y,JX)
        +\langle\overline{\nabla}\psi,(\overline{\nabla}_YJ)X-(\overline{\nabla}_XJ)Y\rangle.
\end{equation}
Here, $\langle\cdot,\cdot\rangle$ is any Riemannian metric on $M$ and $\overline{\nabla}$ is its Levi-Civita connection.
\end{enumerate}
\end{lemma}

\vspace{.1in}

\noindent Now we turn to $J$-holomorphic curves.

\begin{definition}\label{def-holomorphic}
Let $(M^{2n},J)$ be an almost complex manifold and $\Sigma$ be a surface. We call an immersion $F:\Sigma \to (M,J)$ \textbf{$J$-holomorphic} if $J_{F(x)}$ maps $F_{*x}(T_x\Sigma)$ onto itself for any point $x\in\Sigma$.
\end{definition}

\noindent  Stokes' theorem immediately gives the following

\begin{proposition}
If $F:\Sigma\to (M,J)$ be a $J$-holomorphic immersion then $\mathcal{A}$ is constant on $\mathcal{H}$, in particular any $\bar{\omega}_{\rho}\in (\Lambda^2 M)^+\cap [\bar\omega]$ is both a critical point and a stable point for the area functional ${\mathcal A}$.
\end{proposition}

\textbf{Proof:} By the definition of $J$-holomorphic immersion, we can easily see that the almost complex structure $J$ on $M$ can induce an almost complex structure $j$ on $\Sigma$, such that the immersion $F:(\Sigma,j)\to(M,J)$ is $(j,J)$-holomorphic. That if,
\begin{equation}\label{jJ}
    J\circ F_*=F_*\circ j.
\end{equation}
Given any curve $\bar\omega(t)=\bar\omega+d\beta(t)$ which is tamed by $J$, where $\beta(t)$ is a family of smooth 1-forms on $M$, we define the associated Riemannian metric $\bar g(t)$ by (\ref{e2.1}). It suffices to show that
\begin{equation*}
    {\mathcal A}'(t)=0
\end{equation*}
for each $t$. In order to show this, for fixed $t$, at a given point $x$, we take  local coordinates $\{x_1,x_2\}$ on $\Sigma$ such that $\{\partial_{x_1},\partial_{x_2}\}$ is $g(t)$-orthonormal. Then by the above computation, we get that
\begin{equation*}\label{1-1st}
    {\mathcal A}'(t)=\frac{1}{2}\sum_{i=1}^{2}\int_{\Sigma}\bar\omega'(t)\left(F_{x_i},J F_{x_i}\right)d\mu.
\end{equation*}
By the choice of the local frame, it is easy to see that at the given point, $j\partial_{x_1}=\pm\partial_{x_2}$. Without loss of generality, we assume that $j\partial_{x_1}=\partial_{x_2}$, $j\partial_{x_2}=-\partial_{x_1}$. By (\ref{jJ}), we have
\begin{eqnarray*}
\frac{1}{2}\sum_{i=1}^{2}\bar\omega'(t)\left(F_{x_i},J F_{x_i}\right)
   &=& \frac{1}{2}[d\beta'(t)(F_*(\partial_{x_1}),J F_*(\partial_{x_1}))+d\beta'(t)(F_*(\partial_{x_2}),J F_*(\partial_{x_2}))] \\
   &=& \frac{1}{2}[d\beta'(t)(F_*(\partial_{x_1}),F_*j(\partial_{x_1}))+d\beta'(t)(F_*(\partial_{x_2}),F_*j(\partial_{x_2}))] \\
   &=& \frac{1}{2}[(F^*d\beta'(t))(\partial_{x_1},\partial_{x_2})+(F^*d\beta'(t))(\partial_{x_2},-\partial_{x_1})] \\
   &=& (d (F^*\beta'(t)))(\partial_{x_1},\partial_{x_2}).
\end{eqnarray*}
Therefore, at this point,
\begin{equation*}
    \frac{1}{2}\sum_{i=1}^{2}\bar\omega'(t)\left(F_{x_i},J F_{x_i}\right)d\mu=(d (F^*\beta'(t)))(\partial_{x_1},\partial_{x_2})dx_1\wedge dx_2
    =d (F^*\beta'(t)),
\end{equation*}
which is a globally defined exact 2-form on $\Sigma$. As $\Sigma$ is closed, by Stokes' theorem, we see that ${\mathcal A}'(t)=0$. This proves the theorem.
\hfill Q.E.D.

\vspace{.2cm}

\noindent Our interest is in whether (and in which sense) the converse holds.

\vspace{.2cm}

\noindent Let $(M^{2n},\bar{\omega},J,\bar g)$ be a symplectic manifold with symplectic form $\bar\omega$, compatible almost complex structure $J$ and associated Riemannian metric $\bar g$. Recall that the K\"ahler angle $\alpha$ of a surface $\Sigma^2$ in $M$ is defined by (\cite{CW})
\begin{equation*}
    \bar\omega|_{\Sigma}=\cos\alpha d\mu_{\Sigma},
\end{equation*}
where $d\mu_{\Sigma}$ is the induced volume form on $\Sigma$. The following fact is well known:

\begin{proposition}\label{prop2.2}
Let $(M^{2n},\bar{\omega},J,\bar g)$ be a symplectic manifold with compatible almost complex structure $J$. Then $F:\Sigma \to M$ is $J$-holomorphic if and only if $\sin\alpha\equiv0$.
\end{proposition}

\noindent Our main result in this section is as follows:

\begin{theorem}\label{stationary-holomorphic}
Let $(M^{2n},\bar{\omega},J)$ be a compact symplectic manifold with compatible almost complex structure $J$ and $F:\Sigma^2\to M$ be an immersion. If the area functional ${\mathcal A}$ has a critical point in $\mathcal{H}$, then the immersion is $J$-holomorphic.
\end{theorem}

\textbf{Proof:} By definition, there exists a smooth function $\rho$ on $M$, such that $\bar{\omega}_{\rho}(0)=\bar{\omega}_{\rho}=\bar\omega+dd^c\rho\in (\Lambda^2 M)^+\cap [\bar\omega]$ and
\begin{equation*}
    {\mathcal A}'(0)=0
\end{equation*}
for any 
\begin{equation}
\label{omega}
\bar\omega_{\rho}(t)=\bar{\omega}_{\rho}+dd^c\varphi(t)\in (\Lambda^2 M)^+\cap [\bar\omega] \,\, \mbox{with} \,\,  \varphi(0)=0 \,\, .
\end{equation} We will first express the first variation formula (\ref{first}) in terms of K\"ahler angle. Note that in general, $J$ does not need be compatible with $\bar\omega_{\rho}$. We denote $\alpha$ the K\"ahler angle define by $(\bar\omega,J,\bar g)$.

\noindent Fix a point $x\in \Sigma$, it is easy to see that we can choose a $\bar g$-orthonormal frame $\{e_1, e_2, \cdots, e_{2n}\}$ of $T_xM$, such that $\{e_1, e_2\}$ spans the tangent space of $\Sigma$, $\{e_3, \cdots, e_{2n}\}$ spans the normal space of $\Sigma$, and the almost complex structure takes the form
\begin{equation}\label{complex}
J=\left(
  \begin{array}{cccc}
          (J_1)_{4\times 4}     &        0_{4\times (2n-4)}          \\
            0_{(2n-4)\times 4}     &    (J_2)_{(2n-4)\times (2n-4)}        \\
  \end{array}
\right),
\end{equation}
where
\begin{equation}\label{complex1}
J_1=\left(
  \begin{array}{cccc}
          0     &    \cos\alpha    &   \sin\alpha    &        0          \\
    -\cos\alpha &         0        &        0        &   -\sin\alpha     \\
    -\sin\alpha &         0        &        0        &    \cos\alpha     \\
          0     &    \sin\alpha    &   -\cos\alpha   &        0          \\
  \end{array}
\right),
\end{equation}
and $J_2$ satisfies $J_2^2=-Id_{2n-4}$. Suppose $\bar\omega_{\rho}(t)=\bar\omega_{\rho}+dd^c\varphi(t)$, then by (\ref{first}), we have
\begin{equation}\label{1st}
    {\mathcal A}'(0)=\frac{1}{2}\sum_{i,j=1}^{2}\int_{\Sigma}g_{\rho}^{ij}(0)(dd^c\psi)\left(e_i,Je_j\right)d\mu_{\rho},
\end{equation}
where $\psi=\varphi'(0)$, $g_{\rho,ij}=g_{\rho}(e_i,e_j)$ and $(g_{\rho}^{ij})=(g_{\rho,ij})^{-1}$. Plugging (\ref{complex}) into (\ref{1st}) and using (\ref{e2.13}), we finally get that
\begin{eqnarray}\label{1st-new}
 {\mathcal A}'(0)
&=& \frac{1}{2}\int_{\Sigma}g_{\rho}^{11}
               \left\{(\overline{\nabla}^2\psi)(e_{1},e_{1})+\cos\alpha^2(\overline{\nabla}^2\psi)(e_{2},e_{2})
                  +\sin^2\alpha(\overline{\nabla}^2\psi)(e_{3},e_{3})\right.\nonumber\\
   & &  \ \ \ \ \ \ \ \left.   +2\sin\alpha\cos\alpha(\overline{\nabla}^2\psi)(e_{2},e_{3})\right\}
                 d\mu_{\rho} \nonumber\\
   & & +\frac{1}{2}\int_{\Sigma}g_{\rho}^{22}
               \left\{(\overline{\nabla}^2\psi)(e_{2},e_{2})+\cos^2\alpha(\overline{\nabla}^2\psi)(e_{1},e_{1})
                  +\sin^2\alpha(\overline{\nabla}^2\psi)(e_{4},e_{4})\right.\nonumber\\
   & & \ \ \ \ \ \ \ \left.+2\sin\alpha\cos\alpha(\overline{\nabla}^2\psi)(e_{1},e_{4})\right\}
                 d\mu_{\rho} \nonumber\\
   & & +\int_{\Sigma}g_{\rho}^{12}
               \left\{\sin^2\alpha(\overline{\nabla}^2\psi)(e_{1},e_{2})-\sin^2\alpha(\overline{\nabla}^2\psi)(e_{3},e_{4})\right.\nonumber\\
   & & \ \ \ \ \ \ \ \left.  -\sin\alpha\cos\alpha[(\overline{\nabla}^2\psi)(e_{1},e_{3})+(\overline{\nabla}^2\psi)(e_{2},e_{4})]\right\}
                 d\mu_{\rho}\nonumber\\
   & & +\frac{1}{2}\int_{\Sigma}
               \left\{\langle\overline{\nabla}\psi,g_{\rho}^{11}[(\overline{\nabla}_{Je_1}J)e_1-(\overline{\nabla}_{e_1}J)(Je_1)]
                 +g_{\rho}^{22}[(\overline{\nabla}_{Je_2}J)e_2-(\overline{\nabla}_{e_2}J)(Je_2)] \rangle\right\}d\mu_{\rho} \nonumber\\
   & & +\frac{1}{2}\int_{\Sigma}
               \left\{\langle\overline{\nabla}\psi,g_{\rho}^{12}[(\overline{\nabla}_{Je_1}J)e_2-(\overline{\nabla}_{e_2}J)(Je_1)
               +(\overline{\nabla}_{Je_2}J)e_1-(\overline{\nabla}_{e_1}J)(Je_2)] \rangle\right\}d\mu_{\rho}.
\end{eqnarray}
Here, $\langle\cdot,\cdot\rangle=\bar g$ and $\overline{\nabla}$ is its Levi-Civita connection.
By Proposition \ref{prop2.2}, it suffices to show that $\sin\alpha\equiv0$ on $\Sigma$. We prove this by taking special $\psi$ in the first variation formula (\ref{1st-new}). We identify $\Sigma$ with its image in $M$. Denote $d$ the distance function of $M$ from $\Sigma$ with respect to the metric $\bar g$. Namely, for $Q\in M$, $d(Q)=dist_{\bar g}(Q,\Sigma)$. Then it is known that (\cite{Ma}) $\eta=\frac{1}{2}d^2$ is smooth in a neighborhood of $\Sigma$ in $M$.

\vspace{.1in}

\noindent The following result is known for $M=\textbf{R}^{2n}$ (Theorem 3.1 of \cite{AS}), and it is easy to prove using computations in \cite{Ma}.

\vspace{.1in}

\begin{proposition}\label{prop2.6}
Let $\Sigma$ be a $C^{\infty}$ regular submanifold of a $C^{\infty}$ Riemannian manifold $M$, then for any $x_0\in S$, the hessian $Hess(\eta)(x_0)=\frac{1}{2}Hess(d^2)(x_0)$ represents the orthogonal projection on the normal space to $S$ at $x_0$. Namely, for each $X,Y\in T_{x_0}M$, we have
\begin{equation}\label{AB1.3}
    Hess(\eta)(X,Y)(x_0)=\langle X^{\perp},Y^{\perp}\rangle,
\end{equation}
where $T_{x_0}M=T_{x_0}S\oplus N_{x_0}S$ and $X^{\perp}$ is the projection of $X$ onto $N_{x_0}S$.
\end{proposition}

\vspace{.1in}

\noindent Now we can finish the proof of Theorem \ref{stationary-holomorphic}.

\vspace{.1in}

\textbf{Proof of Theorem \ref{stationary-holomorphic} (continued):} We take $\psi$ to be a smooth function on $M$ such that $\psi=\eta$ in a neighborhood of $\Sigma$ in $M$. Then we have that $\overline{\nabla}\psi$ vanishes restricting on $\Sigma$. By the choice of the frame and Proposition \ref{prop2.6}, for $\bar\omega_{\rho}(t)=\bar\omega_{\rho}+dd^c\varphi(t)$ with $\varphi'(0)=\psi$ as above, we have from (\ref{1st-new})
\begin{equation*}
    {\mathcal A}'(0)=\frac{1}{2}\int_{\Sigma}\sin^2\alpha\left(g_{\rho}^{11}+g_{\rho}^{22}\right) d\mu_{\rho}.
\end{equation*}
By our assumption, we must have ${\mathcal A}'(0)=0$ for any $\psi$. In particular, for this special choice of $\psi$, it implies that $\sin\alpha\equiv0$. This proves the theorem.
\hfill Q.E.D.

\vspace{.2in}

\section{Compatible stable point and $J$-holomorphicity}

\vspace{.1in}

\noindent In this section, we will compute the second variation formula for the functional ${\mathcal A}$ (not necessarily at a critical point). Using this formula, we show that if ${\mathcal A}$ has a compatible stable point, then the immersion is $J$-holomorphic.
\vspace{.1in}

\noindent Let $\bar{\omega}_{\rho}(t)$ be a variation of $\bar{\omega}_{\rho}$ as in (\ref{omega}). By (\ref{e2.6}), we have
\begin{equation}\label{e4.1}
    \frac{d^2}{dt^2}|_{t=0}{\mathcal A}(t)=\int_{\Sigma}\frac{d^2}{dt^2}|_{t=0}\nu_{\rho}(t)\sqrt{det(g_{\rho,ij}(0))},
\end{equation}
where $\nu(t)$ is defined by (\ref{nu}). To evaluate $\frac{d^2}{dt^2}|_{t=0}\nu_{\rho}(t)$ at a given point $x$, we  choose the coordinate system $g_{\rho}$-orthonormal at $x$. Thus, we have
\begin{equation}\label{e4.2}
    \frac{d}{dt}\nu_{\rho}(t)=\frac{1}{2}g^{ij}_{\rho}(t)\frac{d}{dt}g_{\rho,ij}(t) \nu_{\rho}(t),
\end{equation}
and
\begin{equation}\label{e4.3}
     \frac{d^2}{dt^2}|_{t=0}\nu_{\rho}(t)=\frac{1}{2}\sum_{i=1}^{2}g_{\rho,ii}''(0)-\frac{1}{2}\sum_{i,j=1}^{2}\left(g_{\rho,ij}'(0)\right)^2
         +\frac{1}{4}\left(\sum_{i=1}^{2}g_{\rho,ii}'(0)\right)^2.
\end{equation}
\noindent By (\ref{e2.9}), we have
\begin{equation}\label{e4.5}
 \sum_{i=1}^{2}g_{\rho,ii}'(0)=\sum_{i=1}^{2}\bar\omega_{\rho}'(0)(e_i,J e_i).
\end{equation}
\noindent By (\ref{metric}), we have
\begin{equation}\label{1-metric}
    g_{\rho,ij}'(t)=\frac{1}{2}\left[\bar\omega_{\rho}'(t)(F_{x_i},J F_{x_j})+\bar\omega'_{\rho}(t)(F_{x_j},J F_{x_i})\right]
\end{equation}
so that
\begin{equation}\label{e4.6}
        g_{\rho,ij}'(0)=\frac{1}{2}\left[\bar\omega_{\rho}'(t)(e_i,J e_j)+\bar\omega_{\rho}'(t)(e_j,J e_i)\right]
\end{equation}
\noindent By (\ref{1-metric}), we have
\begin{equation}\label{e4.7}
    \sum_{i=1}^{2}g_{\rho,ii}''(0)=\sum_{i=1}^{2}\bar\omega_{\rho}''(0)(e_i,J e_i).
\end{equation}

\vspace{.1in}

\noindent Combining (\ref{e4.1}), (\ref{e4.3}), (\ref{e4.5}), (\ref{e4.6}) and (\ref{e4.7}) together, we have
\begin{eqnarray}\label{2-2nd}
 {\mathcal A}''(0)
   &=& \frac{1}{2}\sum_{i=1}^{2}\int_{\Sigma} \bar\omega_{\rho}''(0)(e_i,J e_i)d\mu\nonumber\\
   & & -\frac{1}{8}\sum_{i,j=1}^{2}\int_{\Sigma}
        \left[\bar\omega_{\rho}'(0)(e_i,J e_j)+\bar\omega_{\rho}'(0)(e_j,J e_i)\right]^2d\mu\nonumber\\
   & & +\frac{1}{4}\int_{\Sigma}
         \left\{\sum_{i=1}^2\bar\omega_{\rho}'(0)(e_i,J e_i)\right\}^2d\mu.
\end{eqnarray}
By direct computation, we have
\begin{eqnarray*}
   & & \sum_{i,j=1}^{2}\left[\bar\omega_{\rho}'(0)(e_i,J e_j)+\bar\omega'_{\rho}(0)(e_j,J e_i)\right]^2
        -2\left\{\sum_{i=1}^{2}\left[\bar\omega_{\rho}'(0)(e_i,J e_i)\right]\right\}^2\nonumber\\
   &=& 2\left[\bar\omega_{\rho}'(0)(e_1,J e_2)+\bar\omega_{\rho}'(0)(e_2,J e_1)\right]^2
        +2\left[\bar\omega'_{\rho}(0)(e_1,J e_1)-\bar\omega'_{\rho}(0)(e_2,J e_2)\right]^2.
\end{eqnarray*}
Plugging this into (\ref{2-2nd}) we get
\begin{eqnarray}\label{2-2nd2}
 {\mathcal A}''(0)
   &=& \frac{1}{2}\sum_{i=1}^{2}\int_{\Sigma} \bar\omega_{\rho}''(0)(e_i,J e_i)d\mu\nonumber\\
   & & -\frac{1}{4}\int_{\Sigma}\left[\bar\omega_{\rho}'(0)(e_1,J e_2)+\bar\omega_{\rho}'(0)(e_2,J e_1)\right]^2d\mu \nonumber\\
   & & -\frac{1}{4}\int_{\Sigma}\left[\bar\omega_{\rho}'(0)(e_1,J e_1)-\bar\omega_{\rho}'(0)(e_2,J e_2)\right]^2d\mu.
\end{eqnarray}

\vspace{.1in}

\noindent In particular, when $\bar\omega_{\rho}(t)$ is given by (\ref{omega}), defining
\begin{equation*}
    \frac{\partial\varphi}{\partial t}|_{t=0}=\psi, \ \ \frac{\partial^2\varphi}{\partial t^2}|_{t=0}=\eta,
\end{equation*}
we get
\begin{equation*}
    \frac{\partial\bar{\omega}_{\rho}(t)}{\partial t}|_{t=0}=dd^c\psi, \ \ \frac{\partial^2\bar{\omega}_{\rho}(t)}{\partial t^2}|_{t=0}=dd^c\eta.
\end{equation*}
Then the second variation formula reads
\begin{eqnarray}\label{2nd2}
 {\mathcal A}''(0)
   &=& \frac{1}{2}\sum_{i=1}^{2}\int_{\Sigma} \left[(dd^c\eta)(e_i,J e_i)\right]d\mu_{\rho}\nonumber\\
   & & -\frac{1}{4}\int_{\Sigma}\left[(dd^c\psi)(e_1,J e_2)+(dd^c\psi)(e_2,J e_1)\right]^2d\mu_{\rho} \nonumber\\
   & & -\frac{1}{4}\int_{\Sigma}\left[(dd^c\psi)(e_1,J e_1)-(dd^c\psi)(e_2,J e_2)\right]^2d\mu_{\rho}.
\end{eqnarray}

\vspace{.1in}

\begin{remark}\label{rmk4.2}
When $M$ is a complete noncompact symplectic manifold and $\Sigma$ is a complete submanifold, we can follow the same way to compute the first variation and second variation formulas and give similar definitions as in Definition \ref{stationary} and Definition \ref{A-stable}. In this case, we need the test function for the variations of the target metric $\psi$ and $\eta$ to have compact support on $M$.
\end{remark}

\noindent Our main result in this section is:

\begin{theorem}\label{pseudoholomorphic}
Let $(M^{2n},\bar{\omega},J)$ be a compact symplectic manifold with compatible almost complex structure $J$ and $F:\Sigma^2\to M$ be an immersion. If the area functional ${\mathcal A}$ has a compatible stable point in $[\bar\omega]$, then the immersion is $J$-holomorphic.
\end{theorem}

\textbf{Proof:} By definition, there exists a smooth function $\rho$ on $M$, such that $\bar{\omega}_{\rho}(0)=\bar{\omega}_{\rho}=\bar\omega+dd^c\rho\in (\Lambda^2 M)^+\cap [\bar\omega]$ and
\begin{equation*}
    {\mathcal A}''(0)\geq0
\end{equation*}
for any $\bar\omega_{\rho}(t)=\bar{\omega}_{\rho}+dd^c\varphi(t)\in (\Lambda^2 M)^+\cap [\bar\omega]$ with $\varphi(0)=0$. As in the previous section, we need to express the second variation formula in terms of the K\"ahler angle of $\Sigma$ in $M$. As $J$ is compatible with $\bar\omega_{\rho}$ by our assumption, we can define K\"ahler angle $\alpha_{\rho}$ using $(\bar\omega_{\rho},J,\bar g_{\rho})$.

\noindent At a fixed point $x$ on $\Sigma$, we can take a $\bar g_{\rho}$-orthonormal frame $\{e_1, e_2, \cdots, e_{2n}\}$ of $TM$ such that the almost complex structure takes the form (\ref{complex}) (with $\alpha$ replaced by $\alpha_{\rho}$). By direct computation, we have
\begin{eqnarray*}
  D_1
   &:=& (dd^c\psi)(e_1,J e_2)+(dd^c\psi)(e_2,J e_1) \\
   &=&  \cos\alpha_{\rho}[-(dd^c\psi)(e_1,e_1)+(dd^c\psi)(e_2,e_2)]+\sin\alpha_{\rho}[-(dd^c\psi)(e_1,e_4)+(dd^c\psi)(e_2,e_3)]
\end{eqnarray*}
and
\begin{eqnarray*}
 D_2
   &:=& (dd^c\psi)(e_1,J e_1)-(dd^c\psi)(e_2,J e_2)\nonumber\\
   &=& \cos\alpha_{\rho}\left[(dd^c\psi)(e_1,e_2)+(dd^c\psi)(e_2,e_1)\right]
         +\sin\alpha_{\rho}\left[(dd^c\psi)(e_1,e_3)+(dd^c\psi)(e_2,e_4)\right].
\end{eqnarray*}
By (\ref{e2.13}), we know that $(dd^c\psi)(X,Y)+(dd^c\psi)(Y,X)=0$. Therefore,
\begin{equation}\label{D1}
  D_1 =\sin\alpha_{\rho}[-(dd^c\psi)(e_1,e_4)+(dd^c\psi)(e_2,e_3)]
\end{equation}
and
\begin{eqnarray}\label{D2}
 D_2  &=& \sin\alpha_{\rho}\left[(dd^c\psi)(e_1,e_3)+(dd^c\psi)(e_2,e_4)\right]\nonumber \\
      &=&  \sin\alpha_{\rho}\left\{\sin\alpha_{\rho}\left[(\overline{\nabla}_{\rho}^2\psi)(e_1,e_1)-(\overline{\nabla}_{\rho}^2\psi)(e_2,e_2)
           +(\overline{\nabla}_{\rho}^2\psi)(e_3,e_3)-(\overline{\nabla}_{\rho}^2\psi)(e_4,e_4)\right]\right. \nonumber\\
   & & \left.+2\cos\alpha_{\rho}\left[(\overline{\nabla}_{\rho}^2\psi)(e_2,e_3)-(\overline{\nabla}_{\rho}^2\psi)(e_1,e_4)\right]\right.\nonumber\\
   & & \left.+\langle\overline{\nabla}_{\rho}\psi,((\overline{\nabla}_{\rho})_{e_3}J)e_1-((\overline{\nabla}_{\rho})_{e_1}J)e_3
              +((\overline{\nabla}_{\rho})_{e_4}J)e_2-((\overline{\nabla}_{\rho})_{e_2}J)e_4\rangle_{\rho}\right\}.
\end{eqnarray}
where $\langle\cdot,\cdot\rangle_{\rho}=\bar g_{\rho}$ and $\overline{\nabla}_{\rho}$ is its Levi-Civita connection. Here, we used (\ref{e2.13}) again. We will prove Theorem \ref{pseudoholomorphic} by taking special choices of the test function $\psi$. Indeed, we will take the normal extension of some function on $\Sigma$, which we will recall in the following. (For more details, see, for example, Chapter XIV of S. Lang's book \cite{Lang}.)

\vspace{.1in}

\noindent Let $M$ be a $2n$-dimensional Riemannian manifold and $\Sigma$ be a $p$-dimensional submanifold of $M$ with the induced metric. Locally, we can find a function $r>0$ on $\Sigma$ such that if $N_r\Sigma$ denotes the vectors $w$ with norm $||w||< r(x)$ for $w\in N_{x}\Sigma$, then the exponential map
\begin{equation*}
    \exp: N\Sigma \to M
\end{equation*}
given by
\begin{equation*}
    w \mapsto \exp_x(w) \ \ for \ w\in N_{x}\Sigma
\end{equation*}
gives an isomorphism of $N_r\Sigma$ with an open neighborhood of $\Sigma$ in $M$. Given a function $f$ on $\Sigma$, we may extend $f$ to this tubular neighborhood by making $f$ constant in the normal directions, that is, we define
\begin{equation*}
    f_M(\exp_x(w))=f(x).
\end{equation*}
This extension will be called the \textbf{normal extension} of $f$ to a tubular neighborhood of $\Sigma$.

\noindent In the following, we list some properties of $f_M$ without proof. Some proofs of them and more properties can be found in the book \cite{Lang}.

\begin{lemma}\label{lem3.1}

(a) For vector fields $\xi$, $\eta$ on $\Sigma$, we have on $\Sigma$
\begin{equation*}
    (\overline{\nabla}^2f_M)(\xi,\eta)=(\nabla^2f)(\xi,\eta).
\end{equation*}
(b) Let $\nu$ be a normal vector field on $\Sigma$, then on $\Sigma$
\begin{equation*}
    (\overline{\nabla}^2f_M)(\nu,\nu)=0.
\end{equation*}
(c) Let $\nu$ be a normal vector field on $\Sigma$ and $\xi$ be a tangent vector field on $\Sigma$, then
\begin{equation*}
    (\nu \cdot f_M)(x)=0 \ \ \ \ for \ x\in \Sigma
\end{equation*}
and thus
\begin{equation*}
    \xi(\nu \cdot f_M)(x)=0 \ \ \ \ for \ x\in \Sigma.
\end{equation*}
(d) Let $\xi$ be a tangent vector field on $\Sigma$, then
\begin{equation*}
    (\xi \cdot f_M)(x)=(\xi \cdot f)(x) \ \ \ \ for \ x\in \Sigma.
\end{equation*}
\end{lemma}

\vspace{.1in}

\noindent Now we can prove Theorem \ref{pseudoholomorphic}.

\vspace{.1in}

\textbf{Proof of Theorem \ref{pseudoholomorphic} (continued):} By definition, for any $\bar\omega_{\rho}(t)=\bar{\omega}_{\rho}+dd^c\varphi(t)$ with $t$ small, we have
\begin{equation*}
    {\mathcal A}_{\rho}''(0)\geq0.
\end{equation*}
If we take $\bar\omega_{\rho}(t)=\bar{\omega}_{\rho}+tdd^c\psi$, then $dd^c\eta=\bar\omega_{\rho}''(0)=0$. In this case,
\begin{equation*}
    {\mathcal A}_{\rho}''(0)=-\frac{1}{4}\int_{\Sigma}D_1^2d\mu_{\rho}-\frac{1}{4}\int_{\Sigma}D_2^2d\mu_{\rho}\geq 0.
\end{equation*}
Therefore, we must have
\begin{equation}\label{e4.18}
   D_1(\psi)=D_2(\psi) \equiv 0, \ \ \ on \ \Sigma,
\end{equation}
for any $\psi\in C^{\infty}(M,\textbf{R})$. Fix any $q\in \Sigma$, we will prove that $\sin\alpha(q)=0$ by taking special $\psi$ on $M$. At $q$, we choose an $\bar g_{\rho}$-orthonormal frame $\{e_1, e_2, \cdots, e_{2n}\}$ of $T_qM$ such that $J$ takes the form (\ref{complex}). Set
\begin{eqnarray*}
  A(\psi) &=& (\overline{\nabla}_{\rho}^2\psi)(e_1,e_1)-(\overline{\nabla}_{\rho}^2\psi)(e_2,e_2)+(\overline{\nabla}_{\rho}^2\psi)(e_3,e_3)
              -(\overline{\nabla}_{\rho}^2\psi)(e_4,e_4), \\
  B(\psi) &=& 2\left[(\overline{\nabla}_{\rho}^2\psi)(e_2,e_3)-(\overline{\nabla}_{\rho}^2\psi)(e_1,e_4)\right], \\
  C(\psi) &=& \langle\overline{\nabla}_{\rho}\psi,((\overline{\nabla}_{\rho})_{e_3}J)e_1-((\overline{\nabla}_{\rho})_{e_1}J)e_3
              +((\overline{\nabla}_{\rho})_{e_4}J)e_2-((\overline{\nabla}_{\rho})_{e_2}J)e_4\rangle_{\rho}.
\end{eqnarray*}
Then by (\ref{D2}),
\begin{equation}\label{D22}
    D_2(\psi)=\sin\alpha_{\rho}\left[\sin\alpha_{\rho} A(\psi)+\cos\alpha_{\rho} B(\psi)+C(\psi)\right].
\end{equation}
Taking any $f\in C^{\infty}(\Sigma,\textbf{R})$, we have the normal extension $f_M$ of $f$ over a neighborhood of $\Sigma$ in $M$. Let $\psi\in C^{\infty}(M,\textbf{R})$ such that $\psi=f_M$ in a neighborhood of $\Sigma$. We will compute the restrictions of $A(\psi)$, $B(\psi)$ and $C(\psi)$ to $\Sigma$ using Lemma \ref{lem3.1}. By (a) and (b) of Lemma \ref{lem3.1}, we see that on $\Sigma$
\begin{equation*}
    A(\psi)=(\nabla_{\rho}^2f)(e_1,e_1)-(\nabla_{\rho}^2f)(e_2,e_2).
\end{equation*}
Here, $\nabla_{\rho}$ is the Levi-Civita connection of the induced metric $g_{\rho}$. By (c), (d) of Lemma \ref{lem3.1} and Gauss formula, we have on $\Sigma$
\begin{equation*}
    B(\psi)=2\left[(h^{3}_{21}-h^{4}_{11})e_1(f)+(h^{3}_{22}-h^{4}_{12})e_2(f)\right].
\end{equation*}
and
\begin{eqnarray*}
  C(\psi) &=&  \langle ((\overline{\nabla}_{\rho})_{e_3}J)e_1-((\overline{\nabla}_{\rho})_{e_1}J)e_3
              +((\overline{\nabla}_{\rho})_{e_4}J)e_2-((\overline{\nabla}_{\rho})_{e_2}J)e_4, e_1\rangle e_1(f)\\
   & &  \langle ((\overline{\nabla}_{\rho})_{e_3}J)e_1-((\overline{\nabla}_{\rho})_{e_1}J)e_3
              +((\overline{\nabla}_{\rho})_{e_4}J)e_2-((\overline{\nabla}_{\rho})_{e_2}J)e_4, e_2\rangle e_2(f).
\end{eqnarray*}

\vspace{.1in}

\noindent We define a function $f$ around $q$ on $\Sigma$ by
\begin{equation}\label{function}
    f\left(\tilde{\exp}_q(te_1(q)+se_2(q))\right)=t^2-s^2,
\end{equation}
where the exponential map $\tilde{\exp}$ is defined using the induced metric $g_{\rho}$ on $\Sigma$, and then extend $f$ to be a smooth function on the whole $\Sigma$. By definition, we can easily obtain that
\begin{equation*}
    (e_1f)(q)=\frac{d}{dt}|_{t=0}f\left(\tilde{\exp}_q(te_1(q))\right)=\frac{d}{dt}|_{t=0}t^2=0,
\end{equation*}
\begin{equation*}
    (e_2f)(q)=\frac{d}{ds}|_{s=0}f\left(\tilde{\exp}_q(se_2(q))\right)=\frac{d}{ds}|_{s=0}(-s^2)=0.
\end{equation*}
For the second derivative, we have (Page 344 of \cite{Lang})
\begin{equation*}
    (\nabla^2_{\rho} f)(e_1,e_1)(q)=\frac{d^2}{dt^2}|_{t=0}f\left(\tilde{\exp}_q(te_1(q))\right)=\frac{d^2}{dt^2}|_{t=0}t^2=2,
\end{equation*}
\begin{equation*}
    (\nabla^2_{\rho} f)(e_2,e_2)(q)=\frac{d^2}{ds^2}|_{s=0}f\left(\tilde{\exp}_q(se_1(q))\right)=\frac{d^2}{ds^2}|_{s=0}(-s^2)=-2.
\end{equation*}
Taking values on both sides of (\ref{D22}) at $q$ and using (\ref{e4.18}), we have $\sin\alpha_{\rho}(q)=0.$ As $q$ is arbitrary, we know that $\sin\alpha_{\rho}\equiv 0$ on $\Sigma$. Therefore, the immersion is $J$-holomorphic.
\hfill Q.E.D.

\vspace{.1in}

\noindent The proof of the above theorem relies just on local arguments. In fact the same proof works also in the noncompact
case (see Remark \ref{rmk4.2}):

\begin{theorem}\label{noncompactpseudoholomorphic}
Let $(M^{2n},\bar{\omega},J)$ be a complete noncompact symplectic manifold with compatible almost complex structure $J$ and $F:\Sigma^2\to M$ be an immersion. If the area functional ${\mathcal A}$ has a compatible stable point in $[\bar\omega]$, then the immersion is $J$-holomorphic.
\end{theorem}

\vspace{.2in}

\section{The algebraic case: linear projective stability}

\vspace{.1in}

\noindent Let us now assume that the target manifold is an algebraic manifold that embeds into some complex projective space $\textbf{CP}^N$ holomorphically and isometrically, namely that there is an embedding
\begin{equation*}
    \iota: (M,\bar\omega,J,\bar g) \to (\textbf{CP}^N,\omega_{FS}^N,J_{FS},g_{FS}^N),
\end{equation*}
which is holomorphic, such that
\begin{equation}\label{equa1}
    \iota^*\omega_{FS}^N=\bar\omega, \ \ \iota^*g_{FS}^N=\bar g.
\end{equation}

\vspace{.1in}

\noindent Denote by ${\mathcal H}_N$ and ${\mathcal K}_N$ the space of holomorphic vector fields and Killing vector fields on $\textbf{CP}^N$. Then it is well-known that ${\mathcal H}_N={\mathcal K}_N\oplus J{\mathcal K}_N$. Given any $W\in J{\mathcal K}_N$, it will generate a one parameter family of diffeomorphisms $\Phi_t$ of $\textbf{CP}^N$. It is known that there exists a family of smooth functions $\phi(t)$ on $\textbf{CP}^N$, such that $\tilde\omega(t)=\Phi_t^*\omega_{FS}^N=\omega^N_{FS}+dd^c\phi(t)$. Set $\varphi(t)=\phi(t)\circ \iota$, which is a family of smooth functions on $M$. Set $\dot\varphi=\frac{d}{dt}|_{t=0}\varphi(t)$.

\begin{definition}
Given immersion $F:\Sigma^2\to (M,\bar\omega,J,\bar g)$, we call the area functional \textbf{${\mathcal A}$ has a linearly projectively stable point at $\rho \in  \mathcal{H}$} if $\bar{\omega}_{\rho}$ is projectively induced and
\begin{equation*}
    {\mathcal A}''(0)\geq0
\end{equation*}
for any $\bar\omega_{\rho}(t)=\bar{\omega}_{\rho}+tdd^c\dot\varphi$, where $\varphi(t)$ is defined with $\bar\omega$ replaced by $\bar\omega_{\rho}$ as above.
\end{definition}

\begin{theorem}\label{thm-linearly-stable}
Let $(M,\bar\omega,J,\bar g)$ be an algebraic manifold with all structures induced by the projective space as above and $F:\Sigma^2\to M$ be an immersion. If the area functional has a linearly projectively stable point, then the immersion must be holomorphic with respect to the complex structure $J$.
\end{theorem}

\textbf{Proof:} As $J$ is compatible with any K\"ahler metric in $[\bar\omega]$, without loss of generality, we assume that $\rho\equiv0$ so that $\bar\omega_{\rho}=\bar\omega.$
We denote by  $\alpha$ and $\tilde\alpha$ the K\"ahler angle of $F:\Sigma^2\to (M,\bar\omega,J,\bar g)$ and  $\iota\circ F:\Sigma^2\to (\textbf{CP}^N,\omega_{FS}^N,J_{FS},g_{FS}^N)$, respectively. As the embedding is holomorphic and satisfies (\ref{equa1}), we see that $\alpha=\tilde\alpha$. Set $\psi =\dot\varphi$, then by (\ref{2nd2}), the second variation formula for $\bar\omega(t)=\bar\omega+tdd^c \psi$ is given by
\begin{equation}\label{equa2}
    {\mathcal A}''(0)=-\frac{1}{4}\int_{\Sigma}D_1^2d\mu-\frac{1}{4}\int_{\Sigma} D_2^2d\mu,
\end{equation}
where
\begin{equation*}
    D_1 =\sin\alpha[-(dd^c\psi)(e_1,e_4)+(dd^c\psi)(e_2,e_3)]
\end{equation*}
and
\begin{equation}\label{equa3}
   D_2=\sin\alpha\left[(dd^c\psi)(e_1,e_3)+(dd^c\psi)(e_2,e_4)\right].
\end{equation}
By our assumption, we must have
\begin{equation}\label{equa4}
   D_1(W)=D_2(W)\equiv 0.
\end{equation}
Set $\tilde{\psi}=\dot\phi$ and define
\begin{equation}\label{equa6}
   \hat D_2(W)=\sin\alpha\left[(dd^c\tilde\psi)(\tilde e_1,\tilde e_3)+(dd^c\tilde\psi)(\tilde e_2,\tilde e_4)\right].
\end{equation}
Here, $\{\tilde e_1,\cdots,\tilde e_{2N}\}$ is an orthonormal frame of $\textbf{CP}^{N}$ so that $\tilde e_{\sigma}=\iota_* e_{\sigma}$ for $1\leq \sigma\leq4$. It is easy to see that for this choice of frame, the complex structure $J_{FS}^N$ also takes the form (\ref{complex}). (Recall that $\tilde\alpha=\alpha$.) We can have another expression for $\hat D_2(W)$. From $\tilde\omega(t)=\Phi_t^*\omega_{FS}^N=\omega^N_{FS}+dd^c\phi(t)$, we have $dd^c\tilde\psi=\tilde\omega'(0)=L_{W}\omega^N_{FS}$. By direct computation, we can obtain that
\begin{eqnarray}\label{e6.2}
 \hat D_2 (W)
   &=&  \sin\alpha\left\{\sin\alpha\left[\langle\overline{\nabla}^N_{\tilde e_1}W,\tilde e_1\rangle-\langle\overline{\nabla}^N_{\tilde e_2}W,\tilde e_2\rangle
           +\langle\overline{\nabla}^N_{\tilde e_3}W,\tilde e_3\rangle-\langle\overline{\nabla}^N_{\tilde e_4}W,\tilde e_4\rangle\right]\right. \nonumber\\
   & & \left.+\cos\alpha\left[\langle\overline{\nabla}^N_{\tilde e_2}W,\tilde e_3\rangle+\langle\overline{\nabla}^N_{\tilde e_3}W,\tilde e_2\rangle
          -\langle\overline{\nabla}^N_{\tilde e_1}W,\tilde e_4\rangle-\langle\overline{\nabla}^N_{\tilde e_4}W,\tilde e_1\rangle\right]\right\}.
\end{eqnarray}
Here, $\overline{\nabla}^N$ is the Levi-Civita connection on $(\textbf{CP}^N,g_{FS}^N)$. Suppose $W=JV$ for $V\in{\mathcal K}_N$, then using the fact that
\begin{equation*}
    \langle\overline{\nabla}^N_{\tilde e_i}V,\tilde e_j\rangle+ \langle\overline{\nabla}^N_{\tilde e_j}V,\tilde e_i\rangle=0,
\end{equation*}
we can easily obtain that
\begin{equation}\label{e6.4}
    \hat D_2 (W)=\tilde D_2(JV)=-2\sin\alpha\left[\langle\overline{\nabla}^N_{\tilde e_1}V,\tilde e_3\rangle+\langle\overline{\nabla}^N_{\tilde e_2}V,\tilde e_4\rangle\right].
\end{equation}
In order to proceed further, we need the following key lemma:

\vspace{.1in}

\begin{lemma}\label{lem6.2}
For each point $q\in \Sigma\subset\textbf{CP}^N$, there exists a Killing vector field $V_q\in {\mathcal K}_N$, such that
\begin{equation}\label{e6.5}
    \left[\langle\overline{\nabla}^N_{\tilde e_1}V_q,\tilde e_3\rangle+\langle\overline{\nabla}^N_{\tilde e_2}V_q,\tilde e_4\rangle\right](q)=1,
\end{equation}
and $|\overline{\nabla}^N V_q|(q)\in\left[\frac{\sqrt{2}}{2},\sqrt{2}\right]$. In fact, we can take the value of (\ref{e6.5}) to be any real number by choosing appropriate Killing vector $V_q$.
\end{lemma}

\noindent \textbf{Proof of Lemma \ref{lem6.2}:}
Denote the homogeneous coordinate on $\textbf{CP}^N$ by $[Z_0:\cdots :Z_N]$, and suppose w.l.o.g. that  $q=[1:0:\cdots:0]\in U_0:=\{[Z_0:\cdots:Z_N| Z_0\neq 0]\}$.
Affine coordinates on $U_0$ are given by
\begin{equation*}
    z_j=\frac{Z_j}{Z_0}, \ \ 1\leq j\leq N.
\end{equation*}
Set $z_j=x_j+\sqrt{-1}y_j$. We also assume that the homogeneous coordinate is chosen so that the Fubini-Study metric at $q$ is identity. Namely,
\begin{equation*}
    \langle\frac{\partial}{\partial x_i},\frac{\partial}{\partial x_j}\rangle
    =\langle\frac{\partial}{\partial y_i},\frac{\partial}{\partial y_j}\rangle=\delta_{ij},\ \
    \langle\frac{\partial}{\partial x_i},\frac{\partial}{\partial y_j}\rangle=0.
\end{equation*}

\vspace{.1in}

\noindent Set $\pi: \textbf{C}^{N+1}\backslash \{0\}\to \textbf{CP}^N$ the natural projection, i.e., $\pi(Z_0,\cdots,Z_N)=[Z_0:\cdots: Z_N]$.
In the homogeneous coordinate, the vector field
\begin{equation*}
    \tilde X=\sum_{A,B=0}^{N}a^{AB}Z_A\frac{\partial}{\partial Z_B}
\end{equation*}
is a holomorphic vector field on $\textbf{C}^{N+1}$, where $a^{AB}\in\textbf{C}$. It is known that
\begin{equation*}
    X=\pi_*(\tilde X)
\end{equation*}
is a holomorphic vector field on $\textbf{CP}^N$, and its real part or imaginary part is a Killing vector field on $\textbf{CP}^N$. By direct computation, we obtain that
\begin{eqnarray*}
  X &=& \sum_{j=1}^{N}\left\{-a^{00}z_j-\sum_{i=1}^{N}a^{i0}z_i z_j+a^{0j}+\sum_{i=1}^{N}a^{ij}z_i\right\}\frac{\partial}{\partial z_j} \\
   &:=& \sum_{j=1}^{N}(A_j+\sqrt{-1}B_j)\frac{\partial}{\partial z_j}.
\end{eqnarray*}
Set
\begin{equation*}
    V_1=Re(X)=\frac{1}{2}\sum_{j=1}^{N}\left\{ A_j\frac{\partial}{\partial x_j}+B_j\frac{\partial}{\partial y_j}\right\},\
    V_2=Im(X)=\frac{1}{2}\sum_{j=1}^{N}\left\{ B_j\frac{\partial}{\partial x_j}-A_j\frac{\partial}{\partial y_j}\right\}=-J_{FS}^NV_1.
\end{equation*}
We suppose $V_1$ is a Killing vector field. (The case for $V_2$ is similar.) Suppose $a^{AB}=u^{AB}+\sqrt{-1}v^{AB}$, then it is easy to obtain that
\begin{equation*}
    A_j=-u^{00}x_j+v^{00}y_j-\sum_{i=1}^{N}u^{i0}(x_ix_j-y_i y_j)+\sum_{i=1}^{N}v^{i0}(x_i y_j+x_j y_i)+u^{0j}+\sum_{i=1}^{N}(u^{ij}x_i-v^{ij}y_i),
\end{equation*}
\begin{equation*}
    B_j=-v^{00}x_j-u^{00}y_j-\sum_{i=1}^{N}v^{i0}(x_ix_j-y_i y_j)-\sum_{i=1}^{N}u^{i0}(x_i y_j+x_j y_i) +v^{0j}+\sum_{i=1}^{N}(v^{ij}x_i+u^{ij}y_i).
\end{equation*}

\vspace{.1in}

\noindent As ${\textbf{CP}^N}$ is a symmetric space, we know that for any Killing vector field $V$ and tangent vector field $U$, we have
\begin{equation}\label{e6.6}
    \overline{\nabla}^N_U V=[U,V].
\end{equation}
Suppose
\begin{equation*}
    \tilde e_{\sigma}(q)=\left(\sum_{j=1}^{N}P_{\sigma}^{j}\frac{\partial}{\partial x_j}+\sum_{j=1}^{N}Q_{\sigma}^{j}\frac{\partial}{\partial y_j}\right)(q),
    \ \ \ \ 1\leq\sigma\leq4.
\end{equation*}
Define four vector fields
\begin{equation*}
    \xi_{\sigma}=\sum_{j=1}^{N}P_{\sigma}^{j}\frac{\partial}{\partial x_j}+\sum_{j=1}^{N}Q_{\sigma}^{j}\frac{\partial}{\partial y_j}, \ \ \ \ 1\leq\sigma\leq4.
\end{equation*}
Then
\begin{equation}\label{e6.7}
    (\overline{\nabla}^N_{\tilde e_1}V_1)(q)=(\overline{\nabla}^N_{\xi_1}V_1)(q), \ \
    (\overline{\nabla}^N_{\tilde e_2}V_1)(q)=(\overline{\nabla}^N_{\xi_2}V_1)(q).
\end{equation}
By direct computation, using the definitions of $V_1$, $\xi_1$, $\xi_2$, and recalling that $z_j(q)=0$, we can obtain that
\begin{eqnarray}\label{e6.8}
  [\xi_{\sigma},V_1](q)
   &=& \frac{1}{2}\sum_{j=1}^{N}\left\{P_{\sigma}^j(u^{jj}-u^{00})+\sum_{i\neq j}P_{\sigma}^iu^{ij}
               +Q_{\sigma}^j(v^{00}-v^{jj})-\sum_{i\neq j}Q_{\sigma}^iv^{ij}\right\}\frac{\partial}{\partial x_j} \nonumber\\
   & & +\frac{1}{2}\sum_{j=1}^{N}\left\{-P_{\sigma}^j(v^{00}-v^{jj})+\sum_{i\neq j}P_{\sigma}^iv^{ij}
               +Q_{\sigma}^j(u^{jj}-u^{00})+\sum_{i\neq j}Q_{\sigma}^iu^{ij}\right\}\frac{\partial}{\partial y_j}.
\end{eqnarray}
Written in matrix language, the coordinate of $[\xi_{\sigma},V_1](q)$ in $\{x_i,\cdots,x_N,y_1,\cdots,y_N\}$ is given by
$\frac{1}{2}(P_{\sigma}^1,\cdots,P_{\sigma}^N,Q_{\sigma}^1,\cdots,Q_{\sigma}^N)O$. Here $O$ is a $2N\times 2N$ matrix given by
\begin{equation*}
   O= \left(
  \begin{array}{cc}
    O_1 & O_2 \\
    -O_2 & O_1 \\
  \end{array}
\right),
\end{equation*}
where
\begin{equation*}
   O_1= \left(
  \begin{array}{ccccc}
    u^{11}-u^{00} & u^{12} & \cdot & \cdot & u^{1N} \\
    u^{21} & u^{22}-u^{00} & \cdot & \cdot & u^{2N} \\
    \cdot & \cdot & \cdot & \cdot & \cdot \\
    \cdot & \cdot & \cdot & \cdot & \cdot \\
    u^{N1} & u^{N2} & \cdot & \cdot & u^{NN}-u^{00} \\
  \end{array}
\right)
\end{equation*}
and
\begin{equation*}
   O_2= \left(
  \begin{array}{ccccc}
    v^{11}-v^{00} & v^{12} & \cdot & \cdot & v^{1N} \\
    v^{21} & v^{22}-v^{00} & \cdot & \cdot & v^{2N} \\
    \cdot & \cdot & \cdot & \cdot & \cdot \\
    \cdot & \cdot & \cdot & \cdot & \cdot \\
    v^{N1} & v^{N2} & \cdot & \cdot & v^{NN}-v^{00} \\
  \end{array}
\right).
\end{equation*}
Our goal is to choose an appropriate matrix $O$ (which determines the Killing vector field $V_1$ and in turn the holomorphic vector field $X$) such that $\left[\langle\overline{\nabla}^N_{\tilde e_1}V_1,\tilde e_3\rangle+\langle\overline{\nabla}^N_{\tilde e_2}V_1,\tilde e_4\rangle\right](q)\neq0$. By (\ref{e6.6}), (\ref{e6.7}) and (\ref{e6.8}), we see that
\begin{eqnarray}\label{e6.9}
   & & \left[\langle\overline{\nabla}^N_{\tilde e_1}V_1,\tilde e_3\rangle+\langle\overline{\nabla}^N_{\tilde e_2}V_1,\tilde e_4\rangle\right](q)\nonumber\\
   &=& \frac{1}{2}\left\{\langle(P_{1}^1,\cdots,P_{1}^N,Q_{1}^1,\cdots,Q_{1}^N)O,
             (P_{3}^1,\cdots,P_{3}^N,Q_{3}^1,\cdots,Q_{3}^N)\rangle\right. \nonumber\\
   & & +\left.\langle(P_{2}^1,\cdots,P_{2}^N,Q_{2}^1,\cdots,Q_{2}^N)O,
             (P_{4}^1,\cdots,P_{4}^N,Q_{4}^1,\cdots,Q_{4}^N)\rangle\right\}.
\end{eqnarray}
It is more convenient to write (\ref{e6.9}) in complex notation. Set $\tilde u ^{jj}=u^{jj}-u^{00}$, $\tilde u^{ij}=u^{ij}$ for $i\neq j$. Similarly, set $\tilde v ^{jj}=v^{jj}-v^{00}$, $\tilde v^{ij}=v^{ij}$ for $i\neq j$. Then $O_1=(\tilde u^{ij})$ and $O_2=(\tilde v^{ij})$. By direct computation, using (\ref{e6.8}), we get that
\begin{equation}\label{e6.10}
    [\xi_{\sigma},V_1](q)=Re\left\{\sum_{i,j=1}^{N}(P_{\sigma}^i+\sqrt{-1}Q_{\sigma}^i)
            (\tilde u^{ij}+\sqrt{-1}\tilde v^{ij})\frac{\partial}{\partial z_j}\right\}.
\end{equation}
We also have that
\begin{equation*}
    \tilde e_{\sigma}(q)=2Re\left\{\sum_{j=1}^{N}(P_{\sigma}^j+\sqrt{-1}Q_{\sigma}^j)\frac{\partial}{\partial z_j}\right\}.
\end{equation*}
Therefore,
\begin{eqnarray}\label{e6.11}
   & & \left[\langle\overline{\nabla}^N_{\tilde e_1}V_1,\tilde e_3\rangle+\langle\overline{\nabla}^N_{\tilde e_2}V_1,\tilde e_4\rangle\right](q)\nonumber\\
   &=& 2\left\{\left\langle
       Re\left(\sum_{i,j=1}^{N}(P_{1}^i+\sqrt{-1}Q_{1}^i)(\tilde u^{ij}+\sqrt{-1}\tilde v^{ij})\frac{\partial}{\partial z_j}\right),
            Re\left(\sum_{j=1}^{N}(P_{3}^j+\sqrt{-1}Q_{3}^j)\frac{\partial}{\partial z_j}\right)\right\rangle\right. \nonumber\\
   & & +\left.\left\langle
       Re\left(\sum_{i,j=1}^{N}(P_{2}^i+\sqrt{-1}Q_{2}^i)(\tilde u^{ij}+\sqrt{-1}\tilde v^{ij})\frac{\partial}{\partial z_j}\right),
            Re\left(\sum_{j=1}^{N}(P_{4}^j+\sqrt{-1}Q_{4}^j)\frac{\partial}{\partial z_j}\right)\right\rangle\right\}.
\end{eqnarray}
Set $R_{\sigma}=(R_{\sigma}^1,\cdots,R_{\sigma}^n)=(P_{\sigma}^1+\sqrt{-1}Q_{\sigma}^1,\cdots,P_{\sigma}^n+\sqrt{-1}Q_{\sigma}^n)$ and $\tilde{O}=(\tilde u^{ij}+\sqrt{-1}\tilde v^{ij})$. We will take $\tilde O$ as follows:
\begin{equation}\label{e6.12}
    \tilde O=(R_3^1 \bar R_1,\cdots, R_3^n \bar R_1).
\end{equation}
Then
\begin{equation}\label{e6.13}
    |\tilde O|^2=|R_1|^2(|R_3^1|^2+\cdots |R_3^n|^2)=1,
\end{equation}
and
\begin{equation*}
    R_{\sigma}\tilde O=(\sum_{i=1}^{N}R_{\sigma}^i\bar R _1^i)R_3
    =\left(\sum_{i=1}^{N}(P_1^iP_{\sigma}^i+Q_1^iQ_{\sigma}^i)+\sqrt{-1}\sum_{i=1}^{N}(P_1^iQ_{\sigma}^i-Q_1^iP_{\sigma}^i)\right)R_3.
\end{equation*}
In particular, we have
\begin{equation*}
    R_1\tilde O=R_3
\end{equation*}
and
\begin{equation*}
    R_2\tilde O=\sqrt{-1}\sum_{i=1}^{N}(P_1^iQ_{2}^i-Q_1^iP_{2}^i)R_3.
\end{equation*}
By direct computation and using (\ref{e6.11}), we have
\begin{eqnarray}\label{e6.14}
   & & \left[\langle\overline{\nabla}^N_{\tilde e_1}V_1,\tilde e_3\rangle+\langle\overline{\nabla}^N_{\tilde e_2}V_1,\tilde e_4\rangle\right](q)\nonumber\\
   &=& 2\left\{\left\langle
       Re\left(\sum_{j=1}^{N}R_3^j\frac{\partial}{\partial z_j}\right),
            Re\left(\sum_{j=1}^{N}R_3^j\frac{\partial}{\partial z_j}\right)\right\rangle\right. \nonumber\\
   & & +\left.\left\langle
       Re\left(\sqrt{-1}\sum_{i=1}^{N}(P_1^iQ_{2}^i-Q_1^iP_{2}^i)R_3^j\frac{\partial}{\partial z_j}\right),
            Re\left(\sum_{j=1}^{N}R_4^j\frac{\partial}{\partial z_j}\right)\right\rangle\right\}\nonumber\\
   &=& \frac{1}{2}\left\{1+\sum_{i=1}^{N}(P_1^iQ_{2}^i-Q_1^iP_{2}^i)\sum_{i=1}^{N}(P_3^iQ_{4}^i-Q_3^iP_{4}^i)\right\}.
\end{eqnarray}
Recall that $J_{FS}^N(\frac{\partial}{\partial x_j})=\frac{\partial}{\partial y_j}$, $J_{FS}^N(\frac{\partial}{\partial y_j})=-\frac{\partial}{\partial x_j}$. By the definition of $\tilde e_{\sigma}$, it is easy to see that
\begin{equation*}
    \left[\langle\overline{\nabla}_{\tilde e_1}V_1,\tilde e_3\rangle+\langle\overline{\nabla}_{\tilde e_2}V_1,\tilde e_4\rangle\right](q)
    =\frac{1}{2}(1+\langle J_{FS}^N\tilde e_1,\tilde e_2\rangle\langle J_{FS}^N\tilde e_3,\tilde e_4\rangle).
\end{equation*}
By the choice of complex structure (\ref{complex}), we see that
\begin{equation*}
    \langle J_{FS}^N\tilde e_1,\tilde e_2\rangle=\langle J_{FS}^N\tilde e_3,\tilde e_4\rangle=\cos\alpha(q).
\end{equation*}
Therefore
\begin{equation}\label{e6.15}
     \left[\langle\overline{\nabla}^N_{\tilde e_1}V_1,\tilde e_3\rangle+\langle\overline{\nabla}^N_{\tilde e_2}V_1,\tilde e_4\rangle\right](q)
    =\frac{1}{2}(1+\cos^2\alpha(q))\in \left[\frac{1}{2},1\right].
\end{equation}
For the covariant derivative of $V_1$, we note that, by the definition of $V_1$ and the fact that $z_j(q)=0$, we have
\begin{equation}\label{e6.16}
    |\overline{\nabla}^N V_1|^2(q)=\frac{1}{4}\sum_{j=1}^{N}(|\overline{\nabla}^NA_j|^2+|\overline{\nabla}^NB_j|^2)
    =\frac{1}{2}\left\{\sum_{i,j=1}^N((u^{ij})^2+(v^{ij})^2)+(u^{00})^2+(u^{00})^2\right\}.
\end{equation}
We now choose $O$ so that $u^{00}=u^{0j}=u^{i0}=0$, $v^{00}=v^{0j}=v^{i0}=0$ for $1\leq i,j\leq N$. Then by the definition of $\tilde u^{ij}$ and $\tilde v^{ij}$, we obtain that $\tilde u^{ij}=u^{ij}$, $\tilde v^{ij}=v^{ij}$. On the other hand, as $\tilde O=(\tilde u^{ij}+\sqrt{-1}\tilde v^{ij})$, we have by (\ref{e6.13}) that
\begin{equation*}
    \sum_{i,j=1}^N((\tilde u^{ij})^2+(\tilde v^{ij})^2)= |\tilde O|^2=1.
\end{equation*}
Therefore, by (\ref{e6.16}), we see that
\begin{equation}\label{e6.17}
    |\overline{\nabla}^N V_1|^2(q)=\frac{1}{2}\sum_{i,j=1}^n((u^{ij})^2+(v^{ij})^2)=\frac{1}{2}.
\end{equation}
The proof of the lemma is completed for $V_q=V_1$ by replacing $\tilde O$ by $\frac{2}{1+\cos^2\alpha(q)}\tilde O$. From the proof, we can also see that, by multiplying a constant to $\tilde O$, we can take a Killing vector field $V_q$ such that (\ref{e6.5}) takes any real number.
\hfill Q.E.D.

\vspace{.1in}

\noindent The proof of the theorem is by combining Lemma \ref{lem6.2} with the following observation:

\begin{lemma}\label{lem6.3}
Under our assumption, we have $D_2(W)=\hat D_2(W)$.
\end{lemma}

\noindent \textbf{Proof of Lemma \ref{lem6.3}:} By (\ref{equa3}) and (\ref{equa6}), it suffices to prove that
\begin{equation}\label{equa7}
    (dd^c\psi)(e_1,e_3)+(dd^c\psi)(e_2,e_4)=(dd^c\tilde\psi)(\tilde e_1,\tilde e_3)+(dd^c\tilde\psi)(\tilde e_2,\tilde e_4).
\end{equation}
Here, $\{e_1,e_2\}$ is an any orthonormal frame of $F_*(T_q\Sigma)\subset T_{F(q)}M$, and $\{e_3,e_4\}$ is chosen so that the complex structures $J$ takes the form (\ref{complex}). Furthermore, $\tilde e_{\sigma}=\iota_* e_{\sigma}$ for $1\leq \sigma\leq4$. On the other hand, from $\varphi(t)=\phi(t)\circ \iota$, we see that $\psi=\tilde\psi\circ \iota$. Also note that, as $\iota$ is holomorphic and $dd^c=2\sqrt{-1}\partial\bar\partial$, we have $\iota^*dd^c\tilde\psi=dd^c\iota^*\tilde\psi=dd^c\psi$. Thus,
\begin{eqnarray*}
 (dd^c\tilde\psi)(\tilde e_1,\tilde e_3)+(dd^c\tilde\psi)(\tilde e_2,\tilde e_4)
   &=& (dd^c\tilde\psi)(\iota_*e_1,\iota_*e_3)+(dd^c\tilde\psi)(\iota_*e_2,\iota_*e_4)  \\
   &=& (\iota^* dd^c\tilde\psi)(e_1,e_3)+(\iota^* dd^c\tilde\psi)(e_2,e_4) \\
   &=& (dd^c\psi)(e_1,e_3)+(dd^c\psi)(e_2,e_4).
\end{eqnarray*}
\hfill Q.E.D.

\vspace{.1in}

\noindent Now we can complete the proof of Theorem \ref{thm-linearly-stable}. By Lemma \ref{lem6.2}, for each point $q\in \Sigma\subset\textbf{CP}^n$, there is a Killing vector field $V_q\in{\mathcal K}$, such that (\ref{e6.5}) holds. Combining this with (\ref{e6.4}) and Lemma \ref{lem6.3}, we see that $D_2(JV_q)(q)=\hat D_2(JV_q)(q)=-2\sin\alpha(q)$. By (\ref{equa4}), we must have $\sin\alpha(q)=0$. As $q$ arbitrary, we know that $\sin\alpha\equiv 0$ on $\Sigma$. Therefore, the immersion is holomorphic.
\hfill Q.E.D.

\vspace{.1in}

\section{Symplectic Manifolds with rational symplectic forms}

\vspace{.1in}

\noindent In this section, we extend the results in the previous section to the case that the target manifold $M$ is a symplectic manifold with rational symplectic class. Using the approximately $J_M$-holomorphic embedding of $M$ into some complex projective space, we can define the notion of linearly ${\mathcal A}^k$-stable point and prove that if $p=2$, then the existence of linearly ${\mathcal A}^k$-stable pair point implies $J_M$-holomorphicity.

\noindent Let $(M^{2n},\bar\omega,\bar g,J_M)$ be a compact symplectic manifold with symplectic form $\bar\omega$, compatible almost complex structure $J_M$ and associated Riemannian metric $\bar g$, such that for any $X,Y\in TM$,
\begin{equation}\label{EE2.1}
    \bar g(X,Y)=\bar\omega(X,J_MY).
\end{equation}
Since $\bar\omega$ defines a rational cohomology class, then by a Theorem of Borthwick and Uribe (Theorem 1.1 of \cite{BU}), we known that there exists a sequence of embeddings
\begin{equation}\label{EE2.2}
    \iota_k: M\to (\textbf{CP}^{N_k},\omega_{FS}^k,g_{FS}^k,J_{FS}^k),
\end{equation}
such that, if we put
\begin{equation}\label{EE2.3}
    \bar\omega_k=\iota_k^*\omega_{FS}^k, \ \  \bar g_k=\iota_k^*g_{FS}^k,
\end{equation}
then for $k\geq k_0$
\begin{equation}\label{EE2.4}
    \left|\left|\frac{1}{k}\bar\omega_k-\bar\omega\right|\right|_{C^0}\leq\frac{C_1}{k},
\end{equation}
and
\begin{equation}\label{EE2.5}
    \left|\left|\frac{1}{k}\bar g_k-\bar g\right|\right|_{C^0}\leq\frac{C_2}{k},
\end{equation}
for some constants $C_1$ and $C_2$ and large integer $k_0$.

\noindent Recall the definition of K\"ahler angle:
\begin{equation*}
    \bar\omega|_{\Sigma}=\cos\alpha d\mu_g, \ \ \omega^k_{FS}|_{\Sigma}=\cos\alpha_k d\mu_{g_k}.
\end{equation*}
Here, $\alpha_k$ is the K\"ahler angle of the immersion $\iota_k\circ F:\Sigma^2\to (\textbf{CP}^{N_k},\omega_{FS}^k,g_{FS}^k,J_{FS}^k)$. More precisely, \begin{eqnarray*}
  \bar\omega|_{\Sigma}=F^*\bar\omega, & &  \omega^k_{FS}|_{\Sigma}=(\iota_k\circ F)^*\omega^k_{FS}=F^*\bar\omega_k,\\
   g=F^*\bar g, & &  g_k=(\iota_k\circ F)^*g^k_{FS}=F^*\bar g_k.
\end{eqnarray*}
By (\ref{EE2.5}) and the fact that $d\mu_{\frac{1}{k}g_k}=\frac{1}{k}d\mu_{g_k}$, we see that $\frac{1}{k}d\mu_{g_k}\to d\mu_g$ and combining with (\ref{EE2.4}), we see that
\begin{equation}\label{EE2.6}
    \cos\alpha_k\to \cos\alpha, \ \ \sin\alpha_k\to\sin\alpha \ uniformly \ on \ \Sigma.
\end{equation}

\vspace{.1in}

\noindent Set ${\mathcal K}_k$ the space of Killing vector fields on $\textbf{CP}^{N_k}$. Given any holomorphic vector field $W\in J_{FS}^k{\mathcal K}_k$, let $\Phi_t$ be the one-parameter family of diffeomorphisms generated by $W$. Set $\omega^k(t)=\Phi_t^*\omega_{FS}^k=\omega_{FS}^k+dd_{FS}^c\varphi(t)$ for a family of smooth functions $\varphi(t)$ on $\textbf{CP}^{N_k}$.

\vspace{.1in}

\noindent Note that $\frac{1}{k}\bar\omega_k$ and $\bar\omega$ are in the same cohomology class. Thus, there exists a smooth one form $\gamma_k$ on $M$, such that $\bar\omega=\frac{1}{k}\bar\omega_k+d\gamma_k$. We consider a family of projectively induced symplectic forms on $M$ given by
\begin{equation*}
    \bar\omega(t)=\frac{1}{k}\iota_k^*\omega^k(t)=\frac{1}{k}\iota_k^*\Phi_t^*\omega_{FS}^k
    =\frac{1}{k}\bar\omega_k+d(\frac{1}{k}\iota_k^*d_{FS}^c\varphi(t))  \equiv \bar\omega+d\beta_k(t),
\end{equation*}
where $\beta_k(t)=\frac{1}{k}\iota_k^*d_{FS}^c\varphi(t)-\gamma_k$ is a family of smooth 1-forms on $M$. We can then extend the definitions of the associated metrics by (\ref{e2.1}) and the area functional ${\mathcal A}(t)$ by (\ref{defarea}).

\vspace{.1in}

\begin{definition}\label{symp-k-stable}
Given immersion $F:\Sigma^2\to (M,\bar\omega,J,\bar g)$, we call the area functional \textbf{${\mathcal A}$ has a compatible linearly ${\mathcal A}^k$-stable point at $\rho \in \mathcal{H}$} if  $\bar{\omega}_{\rho}$ is compatible with $J$ and
\begin{equation*}
    {\mathcal A}''(0)\geq0
\end{equation*}
for any $\bar\omega(t)=\bar\omega+td\dot \beta_k$, where $\beta_k(t)$ is defined with $\bar\omega$ replaced by $\bar\omega_{\rho}$ in the above construction.
\end{definition}

\vspace{.1in}

\noindent The main result in this section is

\vspace{.1in}

\begin{theorem}\label{symp-k-J-holomorphic}
Let $(M^{2n},\bar{\omega},J_M,\bar g )$ be a symplectic manifold as above and $F:\Sigma^2\to M$ be an immersion. There exists an integer $K_1$, such that if the area functional has a compatible linearly ${\mathcal A}^k$-stable point for some $k\geq K_1$, then the immersion must be $J_M$-holomorphic.
\end{theorem}

\vspace{.1in}

\textbf{Proof:} As $\bar\omega_{\rho}$ is compatible with $J$ by  assumption we may assume, without loss of generality,  $\rho\equiv0$ so that $\bar\omega_{\rho}=\bar \omega$. Suppose $\dot\beta_k=\frac{\partial \beta_k(t)}{\partial t}|_{t=0}=\theta_k$ so that $\bar\omega(t)=\bar\omega+td\theta_k$. Then by (\ref{2-2nd2}), we see that
\begin{eqnarray}\label{EE2.7}
 {\mathcal A}''(0)
   &=& -\frac{1}{4}\int_{\Sigma}\left[(d\theta_k)(e_1,J_M e_2)+(d\theta_k)(e_2,J_M e_1)\right]^2d\mu \nonumber\\
   & & -\frac{1}{4}\int_{\Sigma}\left[(d\theta_k)(e_1,J_M e_1)-(d\theta_k)(e_2,J_M e_2)\right]^2d\mu.
\end{eqnarray}
Set
\begin{equation*}
    D_1=(d\theta_k)(e_1,J_M e_2)+(d\theta_k)(e_2,J_M e_1), \ \ D_2=(d\theta_k)(e_1,J_M e_1)-(d\theta_k)(e_2,J_M e_2).
\end{equation*}
Then we have
\begin{equation*}
  D_1 =\sin\alpha[-(d\theta_k)(e_1,e_4)+(d\theta_k)(e_2,e_3)]
\end{equation*}
and
\begin{equation}\label{EE2.8}
 D_2  = \sin\alpha\left[(d\theta_k)(e_1,e_3)+(d\theta_k)(e_2,e_4)\right].
\end{equation}
Here, at a fixed point $q$ on $\Sigma$, we take the orthonormal frame $\{e_1,\cdots,e_{2n}\}$ such that the almost complex structure $J_M$ takes the form (\ref{complex}).

\vspace{.1in}

\noindent Set
\begin{equation*}
    \tilde D_2=(\omega^k)'(0)(e_1^k,J_{FS}^k e_1^k)-(\omega^k)'(0)(e_2^k,J_{FS}^k e_2^k),
\end{equation*}
where $\{e_1^k,e_2^k\}$ is an any orthonormal frame of $(\iota_k\circ F)_*(T_q\Sigma)\subset T_{(\iota_k\circ F)(q)}{\textbf{CP}^{N_k}}$. It is easy to see that we can choose an orthonormal frame $\{e_1^k,e_2^k,\cdots, e_{2N_k}^k\}$ of $T_{(\iota_k\circ F)(q)}{\textbf{CP}^{N_k}}$ such that $J_{FS}^{N_k}$ takes the form (\ref{complex}).
Then we have two expressions for $\tilde D _2$. On the one hand, for $W=JV$ with $V\in {\mathcal K}_k$, we have (see (\ref{e6.4}))
\begin{equation}\label{EE2.10}
    \tilde D_2 (W)=\tilde D_2(JV)=-2\sin\alpha_k\left[\langle\overline{\nabla}^k_{e_1^k}V,e_3^k\rangle+\langle\overline{\nabla}^k_{e_2^k}V,e_4^k\rangle\right].
\end{equation}
Here $\overline{\nabla}^k$ is the covariant differential of $(\textbf{CP}^{N_k},g_{FS}^k)$. On the other hand, similar to (\ref{EE2.8}), we have
\begin{equation}\label{EE2.11}
    \tilde D_2 (W) = \sin\alpha_k\left[(dd_{FS}^c\psi)(e_1^k,e_3^k)+(dd_{FS}^c\psi)(e_2^k,e_4^k)\right],
\end{equation}
where $\psi=\dot{\varphi}$.

\vspace{.1in}

\noindent We will prove the theorem by contradiction. Suppose $F:\Sigma^2\to (M^{2n},\bar{\omega},\bar g ,J_M)$ is not $J_M$-holomorphic, then there exists a point $q\in\Sigma\subset M$, such that $\sin\alpha(q)\neq 0$. Without loss of generality, we assume that $\sin\alpha(q)=a>0$. Since $\sin\alpha_k\to\sin\alpha$, we know that $\sin\alpha_k(q)>\frac{a}{2}>0$, for $k\geq N_0$ for some integer $N_0$.

\vspace{.1in}

\noindent Next, we will examine the relation between $D_2(W)$ and $\tilde D_2(W)$. Note that by the definition of $\beta_k$, $d\theta_k=d(\frac{1}{k}\iota_k^*d_{FS}^c\psi)$. Therefore, by (\ref{EE2.8}),
\begin{eqnarray}\label{EE2.12}
 D_2(W)(q)  &=& \sin\alpha\left[\frac{1}{k}(\iota_k^*dd_{FS}^c\psi)(e_1,e_3)+\frac{1}{k}(\iota_k^*dd_{FS}^c\psi)(e_2,e_4)\right] \nonumber\\
   &=&  \sin\alpha\left[(dd_{FS}^c\psi)(\frac{1}{\sqrt k}(\iota_k)_*e_1,\frac{1}{\sqrt k}(\iota_k)_*e_3)
      +(dd_{FS}^c\psi)(\frac{1}{\sqrt k}(\iota_k)_*e_2,\frac{1}{\sqrt k}(\iota_k)_*e_4)\right].
\end{eqnarray}
By the choice of the local frame, we know that $\{e_1,e_2\}$ is an any orthonormal frame of $F_*(T_q\Sigma)\subset T_{F(q)}M$, and $\{e_1^k,e_2^k\}$ is an any orthonormal frame of $(\iota_k\circ F)_*(T_q\Sigma)\subset T_{(\iota_k\circ F)(q)}{\textbf{CP}^{N_k}}$. By (\ref{complex}), we see that at $q$
\begin{equation}\label{EE2.13}
    e_3=\frac{J_Me_1-\cos\alpha(q) e_2}{\sin\alpha(q)}, \ \ e_4=-\frac{J_Me_2+\cos\alpha(q)e_1}{\sin\alpha(q)},
\end{equation}
and
\begin{equation}\label{EE2.14}
    e^k_3=\frac{J_{FS}^ke_1^k-\cos\alpha_k(q) e_2^k}{\sin\alpha_k(q)}, \ \ e_4^k=-\frac{J_{FS}^ke_2^k+\cos\alpha_k(q)e_1^k}{\sin\alpha_k(q)}.
\end{equation}
We now fix $e_1$ and $e_2$ and take
\begin{equation*}
    e_1^k=\frac{\frac{1}{\sqrt k}(\iota_k)_*e_1}{|\frac{1}{\sqrt k}(\iota_k)_*e_1|_{g_{FS}^k}}, \ \
    e_2^k=\frac{\frac{1}{\sqrt k}(\iota_k)_*e_2-\langle \frac{1}{\sqrt k}(\iota_k)_*e_2,e_1^k\rangle e_1^k}{|\frac{1}{\sqrt k}(\iota_k)_*e_2-\langle \frac{1}{\sqrt k}(\iota_k)_*e_2,e_1^k\rangle e_1^k|_{g_{FS}^k}}.
\end{equation*}
Note that
\begin{equation*}
    g_{FS}^k(\frac{1}{\sqrt k}(\iota_k)_*e_i,\frac{1}{\sqrt k}(\iota_k)_*e_j)
    =\frac{1}{k}\bar g_k(e_i,e_j)\to \bar g(e_i,e_j)=\delta_{ij}, \ \ as \ k\to\infty.
\end{equation*}
Therefore, it is easy to see that
\begin{equation}\label{EE2.15}
    |\frac{1}{\sqrt k}(\iota_k)_*e_1-e_1^k|_{g_{FS}^k}\to 0, \ \ |\frac{1}{\sqrt k}(\iota_k)_*e_2-e_2^k|_{g_{FS}^k}\to 0, \  \ as \ k\to\infty.
\end{equation}
On the other hand, by Proposition 4.4 of \cite{BU}, we know that
\begin{equation*}
    |\frac{1}{\sqrt k}(\iota_k)_*J_M e_i-\frac{1}{\sqrt k}J_{FS}^k(\iota_k)_* e_i|\leq \frac{C}{\sqrt k}\to 0, \  \ as \ k\to\infty.
\end{equation*}
Combining all the above together, we see that
\begin{equation}\label{EE2.16}
    |\frac{1}{\sqrt k}(\iota_k)_*e_3-e_3^k|_{g_{FS}^k}\to 0, \ \ |\frac{1}{\sqrt k}(\iota_k)_*e_4-e_4^k|_{g_{FS}^k}\to 0, \  \ as \ k\to\infty.
\end{equation}

\vspace{.1in}

\noindent We can now finish the proof. Comparing (\ref{EE2.10}) and (\ref{EE2.11}) and using Lemma \ref{lem6.2}, we see that we can take a Killing vector field $V_k\in {\mathcal K}_k$, such that for $W_k=J_{FS}^kV_k$, we have
\begin{equation*}
    \left[(dd_{FS}^c\psi_k)(e_1^k,e_3^k)+(dd_{FS}^c\psi_k)(e_2^k,e_4^k)\right](q)=-2,
\end{equation*}
for some smooth function $\psi_k$ on $\textbf{CP}^{N_k}$. On the other hand,
\begin{eqnarray*}
   & & |[(dd_{FS}^c\psi_k)(\frac{1}{\sqrt k}(\iota_k)_*e_1,\frac{1}{\sqrt k}(\iota_k)_*e_3)
   +(dd_{FS}^c\psi_k)(\frac{1}{\sqrt k}(\iota_k)_*e_2,\frac{1}{\sqrt k}(\iota_k)_*e_4)]\\
   & & -[(dd_{FS}^c\psi_k)(e_1^k,e_3^k)+(dd_{FS}^c\psi_k)(e_2^k,e_4^k)]| \\
   &\leq& |(dd_{FS}^c\psi_k)(\frac{1}{\sqrt k}(\iota_k)_*e_1-e_1^k,\frac{1}{\sqrt k}(\iota_k)_*e_3)|
         +|(dd_{FS}^c\psi_k)(e_1^k,\frac{1}{\sqrt k}(\iota_k)_*e_3-e_3^k)| \\
   & & +|(dd_{FS}^c\psi_k)(\frac{1}{\sqrt k}(\iota_k)_*e_2-e_2^k,\frac{1}{\sqrt k}(\iota_k)_*e_4)|
        +|(dd_{FS}^c\psi_k)(e_2^k,\frac{1}{\sqrt k}(\iota_k)_*e_4-e_4^k)| \\
   &\leq& C|dd_{FS}^c\psi_k|_{g_{FS}^k}\epsilon_k
\end{eqnarray*}
for some sequence $\epsilon_k\to0$ and constant $C$ independent of $k$. From $\omega^k(t)=\Phi_t^*\omega_{FS}^k=\omega^k_{FS}+dd_{FS}^c\varphi(t)$, we see that $dd_{FS}^c\psi_k=L_{W_k}\omega_{FS}^k$. Therefore,
\begin{equation*}
    |dd_{FS}^c\psi_k|_{g_{FS}^k}\leq 2|\overline{\nabla}^k W_k|_{g_{FS}^k}=2|\overline{\nabla}^k V_k|_{g_{FS}^k}\leq 2\sqrt{2}.
\end{equation*}
In the last inequality, we used Lemma \ref{lem6.2} again. Hence, we have
\begin{eqnarray*}
   & & \left|(dd_{FS}^c\psi_k)(\frac{1}{\sqrt k}(\iota_k)_*e_1,\frac{1}{\sqrt k}(\iota_k)_*e_3)
       +(dd_{FS}^c\psi_k)(\frac{1}{\sqrt k}(\iota_k)_*e_2,\frac{1}{\sqrt k}(\iota_k)_*e_4)\right|(q)\\
   &\geq & |(dd_{FS}^c\psi_k)(e_1^k,e_3^k)+(dd_{FS}^c\psi_k)(e_2^k,e_4^k)|(q)-C|dd_{FS}^c\psi_k|_{g_{FS}^k}(q)\epsilon_k \\
   &\geq& 2-2\sqrt{2}C\epsilon_k\geq 1
\end{eqnarray*}
for $k\geq K_1\geq N_0$. By (\ref{EE2.12}), we see that for $k\geq K_1$, we have
\begin{equation*}
    |D_2(W_k)(q)|\geq a>0.
\end{equation*}
By (\ref{EE2.7}), we see that for $\beta_k$ associated to such $W_k$, we have
\begin{equation*}
    {\mathcal A}''_{\bar \omega}(d\theta_k)<0.
\end{equation*}
This contradicts our assumption and the proof of the Theorem is completed.
\hfill Q.E.D.

\vspace{.2in}

\section{K\"ahler Manifolds with possibly non rational K\"ahler class}

\vspace{.1in}

\noindent We now assume that $(M,J)$ is an algebraic manifold, that is, a submanifold of some complex projective space. When $[\bar\omega]$ is a rational class and $\bar g$ is the metric induced by the Fubini-Study metric Lawson-Simons (Corollary 9 of \cite{LS}) proved that a submanifold of $M$ is holomorphic if the second variation of the area is nonnegative with respect to holomorphic deformation of $M$ in $\textbf{CP}^N$. In Sections 5, we showed in this case that, existence of linearly projectively stable point also implies holomorphicity. In this section we allow $[\bar\omega]$ to be any real K\"ahler class and $\bar g$ any $J$-induced metric. Take any K\"ahler metric $\bar\omega$ on $M$ with $[\bar\omega]\in H^2(M,\textbf{R})\cap H^{1,1}(M,\textbf{C})$. Let $\bar g$ be the Riemannian metric associated to $\bar\omega$ and $J$.

\vspace{.1in}

\noindent As $(M,J)$ is an algebraic manifold it is easy to see that there exists a sequence of K\"ahler forms $\tau_k$ with $[\tau_k]\in H^2(M,\textbf{Q})\cap H^{1,1}(M,\textbf{C})$, such that
\begin{equation}\label{tau_k}
     ||\tau_k-\bar\omega||_{C^{2}}\leq \varepsilon_k,
\end{equation}
with $\varepsilon_k\to 0$ as $k\to\infty$. Here, the $C^2$ norm is taken with respect to the metric $\bar\omega$.
Since $[\tau_k]$ is rational, there exists, for every $k\in \textbf{N}$, a holomorphic line bundle $(L_k,h_k)\to M$ carrying a hermitian connection $D_k$ of curvature $\frac{\sqrt{-1}}{2\pi}D_k^2=\tau_k$. In particular, $c_1(L_k)=[\tau_k]$. For each positive integer $m>0$, the hermitian metric $h_k$ induces a hermitian metric $h_k^m$ on $L_k^m$. Choose an orthonormal basis $\{S_{k,0}^m, \cdots, S_{k,N_{k,m}}^m\}$ of the space $H^0(M,L_k^m)$ of all holomorphic global sections of $L_k^m$. Here, the inner product on $H^0(M,L_k^m)$ is the natural one induced by the K\"ahler metric $\tau_k$ and the hermitian metric $h_k^m$ on $L^m_k$. By Kodaira embedding theorem, there exists an integer $m_{k,0}$ such that if $m\geq m_{k,0}$, then such a basis induces a holomorphic embedding $\Psi_{k,m}$ of $M$ into $\textbf{CP}^{N_{k,m}}$ given by
\begin{equation}\label{Phik}
    \Psi_{k,m}: M\to \textbf{CP}^{N_{k,m}}, \ \ \Psi_{k,m}(z):= [S_{k,0}^m(z): \cdots : S^m_{k,N_{k,m}}(z)].
\end{equation}
Let $\omega^k_{FS}$ be the standard Fubini-Study metric on $\textbf{CP}^{N_{k,m}}$. Then $\frac{1}{m}\Psi_{k,m}^* \omega^k_{FS}$ is a K\"ahler form on $M$ which lies in the same K\"ahler class as $\tau_k$. We call $\frac{1}{m}\Psi_{k,m}^* \omega^k_{FS}$ the Bergman metric. A famous Theorem proved by Tian (\cite{Tian}) tells us that
\begin{equation}\label{tian}
    \left|\left|\frac{1}{m}\Psi_{k,m}^* \omega^k_{FS}-\tau_k\right|\right|_{C^{2}}\leq \frac{C}{\sqrt{m}}.
\end{equation}
Here the $C^2$ norm is taken with respect to the metric $\tau_k$ and the constant $C$ depends on $\tau_k$. Because of (\ref{tau_k}), we can assume that the constant is uniformly bounded with respect to $k$. Although the Bergman metric $\frac{1}{m}\Psi_{k,m}^* \omega^k_{FS}$ depends on the K\"ahler metric $\tau_k$, the set of Bergman metrics
\begin{equation}\label{bergman}
    {\mathcal P}_{k,m}:= \left\{\frac{1}{m}\Psi_{k,m}^*\sigma^*(\omega^k_{FS})| \sigma \in Aut(\textbf{CP}^{N_{k,m}})\right\},
\end{equation}
is independent of the choice of $\tau_k$ in $[\tau_k]$ and ${\mathcal P}_{k}:=\cup^{\infty}_{m=1}{\mathcal P}_{k,m}$ is dense in $[\tau_k]\cap Ka(M)$ in the $C^2$-topology induced by the one on $\Lambda^2M$. Here, $Ka(M)$ is the space of K\"ahler metrics on $M$. It is known that ${\mathcal P}_{k,m}$ has finite dimension for each $k$ and $m$. Set
\begin{equation}\label{bergman-S}
    {\mathcal Q}_{k}:= \left\{\frac{1}{m(k)}\Psi_{k,m(k)}^*\sigma^*(\omega^k_{FS})| \sigma \in Aut(\textbf{CP}^{N_{k,m(k)}})\right\},
\end{equation}
where $m(k)\geq m_{k,0}$ is a sequence of integers such that $m(k)\to\infty$ as $k\to\infty$.

\vspace{.1in}

\noindent Define
\begin{equation*}
    V_k:=\{\bar\omega\}-\{\tau_k\}+{\mathcal Q}_k=\left\{\bar\omega-\tau_k+\frac{1}{m(k)}\Psi_{k,m(k)}^*\sigma^*(\omega^k_{FS})| \sigma \in Aut(\textbf{CP}^{N_{k,m(k)}})\right\}
\end{equation*}
Then $V_k$ is a finitely dimensional submanifold of $[\bar\omega]$. In particular, for any $\sigma(t)\subset Aut(\textbf{CP}^{N_{k,m(k)}})$, there exists a smooth function $\varphi(t)$ on $M$, such that
\begin{equation}\label{omega(t)}
    \bar\omega(t):=\bar\omega-\tau_k+\frac{1}{m(k)}\Psi_{k,m(k)}^*\sigma(t)^*(\omega^k_{FS})=\bar\omega+2\sqrt{-1}\partial\bar\partial \varphi(t)
    =\bar\omega+dd^c\varphi(t).
\end{equation}
Denote by ${\mathcal H}_{N_{k,m(k)}}$ the space of holomorphic vector fields on $\textbf{CP}^{N_{k,m(k)}}$.

\vspace{.1in}

\begin{definition}
Given an immersion $F:\Sigma^2\to (M,\bar\omega,J,\bar g)$, we say that the area functional \textbf{$\mathcal{A}$} has a $k$-linearly projectively stable point  at $\rho \in \mathcal{H}$ if there exists a smooth function $\rho$ on $M$, such that $\bar{\omega}_{\rho}\in Ka(M)$ and
\begin{equation*}
    {\mathcal A}''(0)\geq0
\end{equation*}
for any $\bar\omega(t)=\bar\omega+tdd^c\dot\varphi$, where $\varphi(t)$ is given with $\sigma(0)=id$ and $\bar\omega$ replaced by $\bar\omega_{\rho}$ in the above construction.
\end{definition}

\noindent The main result in this section is as follows:

\begin{theorem}\label{kahk-holomorphic}
Let $(M,J)$ be an algebraic manifold, $\bar\omega$ be any K\"ahler metric and $F:\Sigma^2\to M$ be an immersion. Then there exists an integer $K_2$, such that if if the area functional has a $k$-linearly projectively stable point  at $\rho \in \mathcal{H}$ for some $k\geq K_2$, then the immersion is $J$-holomorphic.
\end{theorem}

\textbf{Proof:} As $J$ is compatible with any K\"ahler metric in $[\bar\omega]$ we assume, without loss of generality,  that $\rho\equiv0$ so that $\bar\omega_{\rho}=\bar\omega.$ By (\ref{2-2nd2}), for $\bar\omega(t) = \bar\omega + t dd^c \dot\varphi$, the second variation formula is given by
\begin{equation}\label{equation1}
    {\mathcal A}''(0)=-\frac{1}{4}\int_{\Sigma}D_1^2d\mu-\frac{1}{4}\int_{\Sigma} D_2^2d\mu,
\end{equation}
where
\begin{equation*}
    D_1=\bar\omega'(0)(e_1,Je_2)+\bar\omega'(0)(e_2,Je_1),
\end{equation*}
\begin{equation*}
    D_2=\bar\omega'(0)(e_1,Je_1)-\bar\omega'(0)(e_2,Je_2).
\end{equation*}
As before, $\{e_1,e_2\}$ is any orthonormal frame of $F_*(T_q\Sigma)\subset T_{F(q)}M$ with respect to the K\"ahler metric $\bar\omega$. Taking any curve $\sigma(t)$ in $Aut(\textbf{CP}^{N_{k,m(k)}})$, let $W$ be the vector field on $\textbf{CP}^{N_{k,m(k)}}$ generating $\sigma(t)$, then it is known that $W\in {\mathcal H}_{N_{k,m(k)}}$. From (\ref{omega(t)}), we see that
\begin{equation*}
    \bar\omega'(0)=dd^c \dot\varphi=\frac{1}{m(k)}\Psi_{k,m(k)}^*L_W(\omega^k_{FS}).
\end{equation*}
Therefore,
\begin{eqnarray}\label{equation2}
  D_2 &=& (L_W\omega^k_{FS})(\frac{1}{\sqrt{m(k)}}(\Psi_{k,m(k)})_* e_1,\frac{1}{\sqrt{m(k)}}(\Psi_{k,m(k)})_* (Je_1)) \nonumber\\
      & & -(L_W\omega^k_{FS})(\frac{1}{\sqrt{m(k)}}(\Psi_{k,m(k)})_* e_2,\frac{1}{\sqrt{m(k)}}(\Psi_{k,m(k)})_* (Je_2))
\end{eqnarray}
On the other hand, for the second variation formula of $\Psi_{k,m(k)}\circ F:\Sigma\to\textbf{CP}^{N_k}$, we have
(see (\ref{e6.2}))
\begin{equation}\label{equation3}
    \hat D_2 = (L_W\omega^k_{FS})(\tilde e_1,J_{FS}^k\tilde e_1)-(L_W\omega^k_{FS})(\tilde  e_2,J_{FS}^k\tilde e_2)
\end{equation}
By the definition of $\Psi_{k,m(k)}$ , we know that $\Psi_{k,m(k)}$ is holomorphic. From (\ref{tau_k}) and (\ref{tian}), we see that the holomorphic embedding $\Psi_{k,m(k)}:(M,\bar\omega)\to (\textbf{CP}^{N_{k,m(k)}},\omega_{FS}^k)$ satisfies that
\begin{equation*}
    \left|\left|\frac{1}{m(k)}\Psi_{k,m(k)}^* \omega^k_{FS}-\bar\omega\right|\right|_{C^{2}}\to 0
\end{equation*}
and
\begin{equation*}
    \left|\left|\frac{1}{m(k)}\Psi_{k,m(k)}^* g^k_{FS}-\bar g\right|\right|_{C^{2}}\to 0.
\end{equation*}
Then following the same argument as in the proof of Theorem \ref{symp-k-J-holomorphic}, the proof is complete.

\hfill Q.E.D.

\vspace{.1in}

\noindent Recently, Popovici (\cite{Po}) proved that, for a K\"ahler manifold with transcendental K\"ahler class (not necessarily algebraic), one can also have a Kodaira-type approximately holomorphic projective embedding theorem and a Tian-type almost-isometry theorem. Similar to what we did above, we can define another notion of $k$-LP stable point and use the same idea to prove that in fact the existence of $k$-LP stable point implies holomorphicity. It is interesting to underline that even when $M$ is algebraic, this and Theorem \ref{kahk-holomorphic} provide two different approximation arguments.

\vspace{.2in}

\begin{appendix}

\section{}

\vspace{.1in}

\noindent In the main part of this paper, we only consider variations of the symplectic form $\bar{\omega}(t)$ in the same cohomology class. In this appendix, we consider the case the target manifold $M$ is any Riemannian manifold. We will show that in this general case, the concept "${\mathcal A}$-stationary" is not well-posed without restrictions on the type of variations of the metric.

\vspace{.1in}

\noindent Suppose $M$ is an $(n+p)$-dimensional Riemannian manifold and $\Sigma$ is a $p$-dimensional submanifold of $M$ and  a family of immersions
\begin{equation*}
    F(t):\Sigma \to (M,\bar{g}(t)).
\end{equation*}
Then the induced metric on $\Sigma$ is given by
\begin{equation*}
    g(t)=F(t)^*\bar{g}(t).
\end{equation*}
Set
\begin{equation}\label{A1}
    {\mathcal A}(t)\equiv Area(F(t)(\Sigma),g(t))=\int_{\Sigma}d\mu(t).
\end{equation}
We will compute the first variation of ${\mathcal A}$.

\vspace{.1in}

\noindent  Let $F_t$ restricted to $\Sigma$ be the variational vector field and $\{x_i\}_{i=1}^{n}$ be the local coordinates on $\Sigma$. Then in local coordinate,
\begin{equation}\label{A2}
    g_{ij}(t)=\bar{g}(t)\left(\frac{\partial F}{\partial x_i}(t),\frac{\partial F}{\partial x_j}(t)\right).
\end{equation}
Denote by $\overline{\nabla}$ and $\nabla$ the Levi-Civita connections on $M$ and $\Sigma$ respectively.
Set
\begin{equation}\label{A3}
    \nu(t)=\frac{\sqrt{det(g_{ij}(t))}}{\sqrt{det(g_{ij}(0))}}.
\end{equation}
Then $\nu(t)$ is well-defined independent of the choice of coordinate system. Furthermore,
\begin{equation}\label{A.4}
    {\mathcal A}(t)=\int_{\Sigma}\sqrt{det(g_{ij}(t))}=\int_{M}\nu(t)\sqrt{det(g_{ij}(0))},
\end{equation}
and therefore
\begin{equation}\label{A.5}
    \frac{d}{dt}|_{t=0}{\mathcal A}(t)=\int_{\Sigma}\frac{d}{dt}|_{t=0}\nu(t)\sqrt{det(g_{ij}(0))}.
\end{equation}
Suppose
\begin{equation}\label{A.6}
    \frac{\partial\bar{g}(t)}{\partial t}|_{t=0}=h.
\end{equation}
To evaluate $\frac{d}{dt}|_{t=0}\nu(t)$ at a given point $x$, we choose an orthonormal coordinate system. Using this and the fact that $\overline{\nabla}_{F_t}F_{x_i}-\overline{\nabla}_{F_{x_i}}F_{t}=[F_t,F_{x_i}]=0$, we get at $x$:
\begin{eqnarray}\label{A.7}
  \frac{d}{dt}|_{t=0}\nu(t)
   &=& \frac{1}{2}\sum_{i=1}^{p}g_{ii}'(0) = \frac{1}{2}\sum_{i=1}^{p}\frac{d}{dt}|_{t=0}
         \left\{\bar{g}(t)\left(\frac{\partial F}{\partial x_i}(t),\frac{\partial F}{\partial x_i}(t)\right)\right\}\nonumber \\
   &=& \frac{1}{2}\sum_{i=1}^{p}h\left(F_{x_i},F_{x_i}\right)+\bar{g}(\overline{\nabla}_{F_t} F_{x_i},F_{x_i})\nonumber \\
   &=& \frac{1}{2}tr_g(F^*h)
        +\sum_{i=1}^{p}\bar{g}(\overline{\nabla}_{F_{x_i}}F_t^T,F_{x_i})
        +\sum_{i=1}^{p}\bar{g}(\overline{\nabla}_{F_{x_i}}F_t^N,F_{x_i})\nonumber \\
   &=& \frac{1}{2}tr_g(F^*h)+div_{\Sigma}F_t^T-\bar{g}(F_t,\textbf{H}),
\end{eqnarray}
Therefore,
\begin{equation}\label{A.8}
    \delta{\mathcal A}(F_t,h)=\frac{1}{2}\int_{\Sigma}tr_g(F^*h) d\mu-\int_{\Sigma}\bar{g}(F_t,\textbf{H})d\mu.
\end{equation}
We cannot expect that ${\mathcal A}'(0)=0$ for any $F_t$ and any $h$. In fact, we have

\begin{proposition}\label{propA1}
Given $(F,\bar{g})$, for any variational vector field $F_t$, there exist $h_1$ and $h_2$, such that $\delta{\mathcal A}(F_t,h_1)>0$ and $\delta{\mathcal A}(F_t,h_2)<0$.
\end{proposition}

\textbf{Proof:} As $F$ and $\bar{g}$ is fixed, given any $F_t$, suppose
\begin{equation*}
    \left|\bar{g}(F_t,\textbf{H})\right| \leq C_0.
\end{equation*}
We take $h_1=2(C_0+1)\bar{g}$ and $h_2=-2(C_0+1)\bar{g}$, then $F^* h_1=2(C_0+1)g$ and $F^* h_2=-2(C_0+1)g$. Therefore,
\begin{equation*}
    \delta{\mathcal A}(F_t,h_1)\geq \left[(p-1)C_0+p\right]Area(\Sigma) \geq pArea(\Sigma)>0
\end{equation*}
and
\begin{equation*}
    \delta{\mathcal A}(F_t,h_2)\leq \left[-(p-1)C_0-p\right]Area(\Sigma)\leq -pArea(\Sigma) <0.
\end{equation*}
\hfill Q.E.D.

\end{appendix}

\end{document}